\input amstex\documentstyle{amsppt}  
\pagewidth{12.5cm}\pageheight{19cm}\magnification\magstep1
\topmatter
\title On the generalized Springer correspondence\endtitle
\author G.Lusztig\endauthor
\thanks{Supported in part by the National Science Foundation}\endthanks
\address{Department of Mathematics, M.I.T., Cambridge, MA 02139}\endaddress
\endtopmatter   
\document
\define\spn{\text{\rm span}}
\define\Irr{\text{\rm Irr}}

\define\dcs{\dot{\cs}}

\define\bcs{\bar{\cs}}

\define\bu{\bar u}

\define\frl{\forall}

\define\si{\sim}

\define\sqc{\sqcup}

\define\qua{\quad}

\define\hP{\hat P}
\define\hS{\hat S}
\define\hV{\hat V}

\define\hcs{\hat{\cs}}

\define\bC{\bar C}
\define\bH{\bar H}
\define\bQ{\bar Q}
\define\bY{\bar Y}

\define\bYY{\ov{\YY}}

\define\spa{\spadesuit}
\define\part{\partial}
\define\emp{\emptyset}
\define\imp{\implies}

\define\n{\notin}

\define\m{\mapsto}
\define\do{\dots}

\define\lra{\leftrightarrow}

\define\sub{\subset}    

\define\T{\times}
\define\ti{\tilde}
\define\nl{\newline}
\redefine\i{^{-1}}

\define\un{\underline}
\define\ov{\overline}
\define\ot{\otimes}
\define\bbq{\bar{\QQ}_l}

\define\Ad{\text{\rm Ad}}

\define\End{\text{\rm End}}

\define\tr{\text{\rm tr}}

\define\di{\diamond}

\define\a{\alpha}
\redefine\b{\beta}
\redefine\c{\chi}
\define\g{\gamma}
\redefine\d{\delta}
\define\e{\epsilon}

\redefine\o{\omega}
\define\p{\pi}
\define\ph{\phi}
\define\ps{\psi}
\define\r{\rho}
\define\s{\sigma}

\define\th{\theta}
\define\k{\kappa}

\define\z{\zeta}
\define\x{\xi}

\define\vt{\vartheta}

\redefine\G{\Gamma}

\define\boc{\bold c}

\define\kk{\bold k}

\define\JJ{\bold J}

\define\NN{\bold N}

\define\PP{\bold P}
\define\QQ{\bold Q}

\define\SS{\bold S}
\define\TT{\bold T}

\define\ZZ{\bold Z}
\define\XX{\bold X}
\define\YY{\bold Y}

\define\cb{\Cal B}

\define\cd{\Cal D}
\define\ce{\Cal E}

\define\ch{\Cal H}
\define\ci{\Cal I}

\define\cl{\Cal L}
\define\cm{\Cal M}

\define\co{\Cal O}
\define\cp{\Cal P}

\define\cs{\Cal S}

\define\cw{\Cal W}
\define\cz{\Cal Z}

\define\fA{\frak A}

\define\fS{\frak S}

\define\fU{\frak U}
\define\fV{\frak V}
\define\fZ{\frak Z}

\define\tj{\ti j}

\define\tq{\ti q}

\define\tx{\ti x}

\define\tM{\ti M}

\define\tP{\ti P}  
\define\tQ{\ti Q}

\define\tS{\ti S}

\define\tV{\ti V}

\define\tY{\ti Y}

\define\sha{\sharp}

\define\bP{\bar P}

\define\DL{DL}
\define\HO{Ho}
\define\ICC{L1}
\define\CSV{L2}
\define\LS{LS}
\define\SPA{S1}
\define\SPAA{S2}
\define\SPR{Sp}
\NoBlackBoxes

\head Introduction\endhead
\subhead 0.1\endsubhead
Let $\kk$ be an algebraically closed field of characteristic $p\ge0$.
Let $G$ be a connected reductive group over $\kk$. 
Let $\fA$ be the set of all pairs $(\boc,\ce)$ where $\boc$ is a unipotent class in $G$ and $\ce$ is an
irreducible $G$-equivariant local system on $\boc$ defined up to isomorphism;
let $\ov{\fA}$ be the set of triples $(P,C,\cs_0)$ (up to $G$-conjugacy) where $P$ is a parabolic
subgroup of $G$, $C$ is a unipotent class of the reductive quotient $\bP$ of $P$ and $\cs_0$ is 
an irreducible $\bP$-equivariant cuspidal local system on $C$.
According to \cite{\ICC, 6.5}, there is a canonical surjective map $\fA@>>>\ov{\fA}$ (whose fibres are 
called blocks) such that the block corresponding to $(P,C,\cs_0)\in\ov{\fA}$ is in natural bijection 
(``generalized Springer correspondence'') with the set $\Irr\cw$ of irreducible representations (up to 
isomorphism) of the finite group $\cw:=NL/L$ where $NL$ is the normalizer in $G$ of a Levi subgroup $L$
of $P$. (One can show that $\cw$ is naturally a Weyl group depending on the block.)
The case considered originally by Springer, see \cite{\SPR}, with some restrictions on $p$, involves
the block corresponding to $(B,\{1\},\bbq)$ where $B$ is a Borel subgroup of $G$.
The problem of determining explicitly the generalized Springer correspondence can be reduced 
to the case where $G$ is almost simple, simply connected. For such $G$ the explicit bijection 
was determined in \cite{\ICC}, \cite{\SPAA} and the references there, for any block except for 

(a) two blocks for $G$ of type $E_6$ with $p\ne3$, with $\cw$ of type $G_2$  and 

(b) two blocks for $G$ of type $E_8$ with $p=3$, with $\cw$ of type $G_2$;
\nl
for these blocks the method of \cite{\SPAA} (based mostly on the restriction theorem 
\cite{\ICC, 8.3}) had the following gap: it gave the explicit bijection only up to composition with 
a permutation of $\Irr\cw$ which interchanges the two $2$-dimensional irreducible 
representations of $\cw$ and keeps fixed all the other irreducible representations.

\subhead 0.2\endsubhead
Let $(\boc,\ce)\in\fA$ and let $(\boc',\ce')\in\fA$ be such that $\boc'$ is contained in the closure
$\bar\boc$ of $\boc$. 
For an integer $k$ let $m^k_{(\boc',\ce'),(\boc,\ce)}$ be the multiplicity of $\ce'$ in
the local system on $\boc'$ obtained by restricting to $\boc'$ the $k$-th cohomology sheaf of the
intersection cohomology complex $IC(\bar\boc,\ce)$. Let
$$m_{(\boc',\ce'),(\boc,\ce)}=\sum_{k\ge0}m^k_{(\boc',\ce'),(\boc,\ce)}q^{k/2}\in\NN[q^{1/2}]$$
where $q^{1/2}$ is an indeterminate. In \cite{\CSV,24.8} it is shown (assuming that $p$ is not a bad
prime for $G$) that $m_{(\boc',\ce'),(\boc,\ce)}\in\NN[q]$,
that $m_{(\boc',\ce'),(\boc,\ce)}$ is $0$ if $(\boc',\ce'),(\boc,\ce)$ are not in the same block and
that $m_{(\boc',\ce'),(\boc,\ce)}$ is explicitly computable in terms of the generalized Springer
correspondence if $(\boc',\ce'),(\boc,\ce)$ are in the same block.

\subhead 0.3\endsubhead
In this subsection we consider a block of $G$ as in 0.1(a),(b) where $p$  not a bad prime for 
$G$; thus we must be in case 0.1(a) and $p\n\{2,3\}$. An attempt to close the gap in this case was
made in \cite{\CSV, 24.10}, based on computing the polynomials
$m_{(\boc',\ce'),(\boc,\ce)}$ (see 0.2) for $(\boc',\ce'),(\boc,\ce)$ in this block.
Unfortunately, the attempt in \cite{\CSV} contained a calculation error. (I thank Frank
L\"ubeck for pointing this out to me.) As a result, the gap remained. In this paper we refine the 
analysis in \cite{\CSV, 24.10} and close the gap for the blocks in 0.1(a) with $p\n\{2,3\}$ (see 
Theorem 5.5). We now describe our strategy. Using the algorithm of \cite{\CSV, 24.10} (based on
the two possible scenarios corresponding to the
two possible inputs for the generalized Springer correspondence) we compute the polynomials
$m_{(\boc',\ce'),(\boc,\ce)}$ for $(\boc',\ce'),(\boc,\ce)$ in the block. We get two sets of
polynomials one in each of the two scenarios. In fact the algorithm gives only rational functions in
$q$ which, by a miracle,  turn out to be in $\NN[q]$ in both scenarios. So by this computation one
cannot rule out one of the scenarios and a further argument is needed. Let $u\in G$ be a unipotent
element of $G$ such that the Springer fibre at $u$ is $3$-dimensional (such $u$ is unique up to
conjugacy) and let $\XX_u$ be the ``generalized Springer fibre'' at $u$ (attached to the block).
It is known from \cite{\ICC} that $\cw$ acts naturally on $H^j_c(\XX_u,?)$ where ? is a suitable
local system on $\XX_u$ attached to the block. By a known result, from the knowledge of the polynomials 
$m_{(\boc',\ce'),(\boc,\ce)}$ one can extract the trace of a simple reflection of $\cw$ on 
$\sum_j(-1)^jH^j_c(\XX_u,?)$. It turns out that this trace is $1$ in one scenario and $-1$ in the
other scenario. We show that this trace is equal to the Euler characteristic of a certain explicit open
subvariety of the Springer fibre at $u$. We then try to compute from first principles this Euler
characteristic; the computation occupies most of the paper. The computation does not give the exact
value of the Euler characteristic; it only shows that it is one of the numbers $0,1,2$. Since we already
know that it is $1$ or $-1$ we
deduce that it is $1$; moreover, this determines which scenario is real and which one is not.

\subhead 0.4\endsubhead
Note that Theorem 5.5 completes the explicit determination of the generalized Springer correspondence
for any $G$ and any block, assuming that $p$ is not a bad prime for $G$. 

In the case where $p$ is a bad prime for $G$, a gap remains, and in 6.2 we state a conjecture about it.

\subhead 0.5\endsubhead
{\it Notation.} All algebraic varieties are assumed to be over $\kk$ and all algebraic groups 
are assumed to be affine.

For any connected algebraic group $H$ let $U_H$ be the unipotent radical of $H$ and let
$\bH=H/U_H$, a connected reductive group; let $\p_H:H@>>>\bH$ be the obvious homomophism.
Let $\cz_{\bH}^0$ be the connected centre of $\bH$.
Let $\cb^H$ (resp. $\cb^{\bH}$) be the variety of Borel subgroups of $H$ (resp. $\bH$). Note 
that $\b\m\p_H\i(\b)$ is an isomorphism $\cb^{\bH}@>\si>>\cb^H$. Now $H$ (resp. $\bH$) acts by
simultaneous conjugation on $\cb^H\T\cb^H$ (resp. $\cb^{\bH}\T\cb^{\bH}$); the orbits of this 
action are naturally parametrized by the Weyl group $W^H$ (resp. $W^{\bH}$) of $H$ (resp. $\bH$).
The identification $\cb^{\bH}\T\cb^{\bH}\lra\cb^H\T\cb^H$ via $\p_H\T\p_H$, induces an 
identification $W^H=W^{\bH}$. 

If $H$ is a connected reductive group we denote by $H_{ad}$ the adjoint group of $H$. In this 
paper $G$ is a fixed connected reductive group. We write $\cb=\cb^G$, $W=W^G$. Let $\co_w$ be the 
$G$-orbit on $\cb\T\cb$ indexed by $w\in W$. The simple reflections in $W$ are denoted as
$\{s_i;i\in I\}$ where $I$ is an indexing set. For any $J\sub I$ let $W_J$ be the subgroup of $W$ 
generated by $\{s_i;i\in J\}$. Let $\e_J:W_J@>>>\{\pm1\}$ be the homomorphism given by 
$\e_J(s_i)=-1$ for all $i\in J$.
The set $\cp$ of parabolic subgroups of $G$ is naturally partitioned
as $\cp=\sqc_{J\sub I}\cp_J$ where for $J\sub I$, $\cp_J$ consists of parabolic subgroups $P$ with
the following property: for $w\in W$, we have $(B,B')\in\co_w$ for some some $B,B'$ in 
$\cb$, $B\sub P$, $B'\sub P$ if and only if $w\in W_J$. Note that $\cb=\cp_\emp$.

Now let $J\sub I$ and let $P\in\cp_J$. Then $\cb^P$ is equal to the closed subvariety
$\{B\in\cb;B\sub P\}$ of $\cb$ and the obvious map $W^P@>>>W$ identifies $W^P$ with the subgroup 
$W_J$ of $W$. The set of simple reflections of $W^P=W^{\bP}$ becomes the set $\{s_i;i\in J\}$.

For $i\in I$ we write $\cp_i$ instead of $\cp_{\{i\}}$ and $\cp^i$ instead of $\cp_{I-\{i\}}$, a 
class of maximal parabolic subgroups. For any $B\in\cb$ there is a unique $P\in\cp^i$ such that 
$B\sub P$; we set $P=B(i)$. A subvariety of $\cb$ is said to be an $i$-line if it is of the form 
$\cb^P$ for some $P\in\cp_i$ (which is necessarily unique). If $P,P'$ are parabolic subgroups of 
$G$ we write $P\spa P'$ whenever $P,P'$ contain a common Borel subgroup.

Let $l$ be a fixed prime number such that $l\ne p$. All local systems are assumed to be
$\bbq$-local systems.

\head Contents\endhead

1. Preliminaries.

2. A trace computation.

3. Computations in certain groups of semisimple rank $\le5$.

4. Euler characteristic computations.

5. The main result.

6. Final comments.

\head 1. Preliminaries\endhead
\subhead 1.1\endsubhead
We fix a unipotent element $u\in G$; let $\cb_u=\{B\in\cb;u\in B\}$. This is a nonempty 
subvariety of $\cb$. Let $\un\cb_u$ be the set of irreducible components of $\cb_u$. According to 
Spaltenstein \cite{\SPA}, for $X\in\un\cb_u$, $d_u:=\dim X$ depends only on $u$, not on $X$. For 
any $X\in\un\cb_u$ let $J_X$ be the set of all $i\in I$ such that $X$ is a union of $i$-lines. We
show:

(a) {\it Let $J=J_X$, let $B\in X$ and let $P$ be the unique subgroup in $\cp_J$ such that $B\sub P$. Then
$\cb^P\sub X$.}
\nl
Let $B'\in X$. We can find a sequence $B=B_0,B_1,\do,B_t=B'$ in $\cb$ such that for 
$k=0,1,\do,t-1$ we have $(B_k,B_{k+1})\in\co_{s_{i_k}}$ where $i_k\in J$. We show by induction on 
$k$ that $B_k\in X$. For $k=0$ this holds by assumption. Assume now that $k\ge1$. By the induction
hypothesis we have $B_{k-1}\in X$. Since $i_k\in J_X$, the $i_k$-line containing $B_{k-1}$ is 
contained in $X$. Since $B_k$ is contained in this $i_k$-line we have $B_k\in X$. This completes 
the induction. We see that $B'\in X$. This proves (a).

We show:

(b) {\it In the setup of (a), assume in addition that $d_u$ is equal to $\nu_J$, the number of 
reflections in $W_J$. Then $\cb^P=X$.}
\nl
Note that $\cb^P$ is a closed irreducible subvariety of dimension $\nu_J$ of $X$ and $X$ is 
irreducible of dimension $d_u$. The result follows.

\subhead 1.2\endsubhead
In this subsection we assume that $G_{ad}=\prod_{i\in I]}^mG_i$ where $G_i\cong PGL_2(\kk)$ for
all $i\in I$. Assume that $B\in\cb_u$ is such that for any $i\in I$, the $i$-line through $B$ is 
contained in $\cb_u$. We show that $u=1$. We write $u=(u_i)$ where $u_i\in G_i$ for all $i$. Let 
$K=\{i\in I;u_i=1\}$. Then $\cb_u$ is a product of copies of $\PP^1$ indexed by $K$. Moreover, 
$\cb_u$ contains an $i$-line if and only if $i\in K$. Our assumption implies that $K=I$ hence 
$u=1$, as asserted.

\subhead 1.3\endsubhead
We shall need the following variant of 1.1(a),(b).

(a) {\it Let $J$ be a subset of $I$ such that $s_js_{j'}=s_{j'}s_j$ for all $j,j'$ in $J$. Let 
$B\in\cb_u$ be such that for any $j\in J$ the $j$-line containing $B$ is contained in $\cb_u$. Let
$P$ be the unique subgroup in $\cp_J$ such that $B\sub P$. Assume that $\sha(J)=d_u$. Then 
$X=\cb^P$ is an irreducible component of $\cb_u$ and $J_X=J$.}
\nl
We apply 1.2 with $G,u$ replaced by $\bP,\p_P(u)$. Note that $\bP_{ad}$ is a product of copies of 
$PGL_2(\kk)$. We see that $\p_P(u)=1$ that is, $u\in U_P$. Then any Borel subgroup of $P$ contains
$u$. Now (a) follows.

\subhead 1.4\endsubhead
Let $[\un\cb_u]$ be the $\bbq$-vector space with basis $\{X;X\in\un\cb_u\}$. Let $\r_u$ be the Springer 
representation of $W$ on the vector space $[\un\cb_u]$. The following property of $\r_u$ appeared in a 
letter of the author to Springer (March 1978), see also \cite{\HO}:

(a) {\it Let $i\in I$ and let $X\in\un\cb_u$. Then

$s_iX=-X$ if $i\in J_X$;

$s_iX-X\in\sum_{X'\in\un\cb_u;i\in J_{X'}}\ZZ X'$ if $i\n J_X$.}
\nl
For $J\sub I$ let $(\un\cb_u)'_J=\{X\in\un\cb_u;J\sub J_X\}$  and let 
$[(\un\cb_u)'_J]$ be the subspace of $[\un\cb_u]$ spanned by $(\un\cb_u)'_J$.
We have $(\un\cb_u)'_J=\cap_{i\in J}(\un\cb_u)'_{\{i\}}$ and
$[(\un\cb_u)'_J]=\cap_{i\in J}[(\un\cb_u)'_{\{i\}}]$.
From (a) we see that for any $i\in J$, we have $[(\un\cb_u)'_{\{i\}}]=\{v\in[\un\cb_u];s_iv=-v\}$.
Hence 
$[(\un\cb_u)'_J]=\cap_{i\in J}\{v\in[\un\cb_u];s_iv=-v\}=\{v\in[\un\cb_u];wv=\e_J(w)v \frl w\in W_Jv\}$.
Thus, $\sha(\un\cb_u)'_J$ is 
equal to the multiplicity $(\e_J:\r_u|_{W_J})$ of $\e_J$ in the $W_J$-module 
$\r_u|_{W_J}$.
Let $(\un\cb_u)_J=\{X\in\un\cb_u;J=J_X\}$
We have $\sha\un\cb_u)'_J=\sum_{J';J\sub J'\sub I}\sha(\un\cb_u)_{J'}$.
Hence $\sha\un\cb_u)_J=\sum_{J';J\sub J'\sub I}(-1)^{|J'|-|J|}\sha\un\cb_u)'_{J'}$ so that

(b) $\sha(\un\cb_u)_J=\sum_{J';J\sub J'\sub I}(-1)^{|J'|-|J|}(\e_{J'}:\r_u|_{W_{J'}})$.

\subhead 1.5\endsubhead
In this subsection we assume that we are given $J\sub I$, $P\in\cp_J$ and $\bu\in\bP$ such that
$u\in P$, $\bu=\p_P(u)$ and such that, if $C$ is the conjugacy class of $u$ in $G$ then

(a) $C\cap\p_P\i(\bu)$ is open dense in $\p_P\i(\bu)$.
\nl
Then $C$ is induced (in the sense of \cite{\LS}) by the $\bP$-conjugacy class of $\bu$.

Let $\cb^P,\cb^{\bP}$ be as in 0.4 and let $\cb^{\bP}_{\bu}=\{\b\in\cb^{\bP};\bu\in\b\}$. Let 
$d_{\bu}=\dim\cb^{\bP}_{\bu}$.  From \cite{\LS} it is known that

(b) $d_u=d_{\bu}$.
\nl
For any irreducible component $\x$ of $\cb^{\bP}_{\bu}$, the subset $J_\x$ of $J$ is defined as
in 1.1 (replacing $G,I,u$ by $\bP,J,\bu$). Let $\ti\x$ be the image of $\x$ under 
$\b\m\p_P\i(\b)$; note that $\ti\x$ is a closed irreducible subvariety of $\cb_u$ of dimension $d_{\bu}$ hence by (b) is an irreducible component of $\cb_u$. From the definitions we see that 

(c) $J_\x=J_{\ti\x}$ where $J_{\ti\x}$ is defined as in 1.1.

\subhead 1.6\endsubhead
Let $i\in I$. Let $\cb_{u,i}$ be the set of all $B\in X$ such that the $i$-line though $B$ is 
contained in $\cb_u$ or equivalently such that if $P_i\in\cp_i$ is defined by $B\sub P_i$ then
$u\in U_{P_i}$.

For $X\in\un\cb_u$ let $X_i=X\cap\cb_{u,i}$.

\head 2. A trace computation\endhead
\subhead 2.1\endsubhead
We fix $J\sub I$ such that $J\ne I$ and $P_J\in\cp_J$. We write $\p_J$ instead of $\p_{P_J}$.
We fix a unipotent conjugacy class $C$ of $\bP_J$ and let $\cs_0$ be an irreducible cuspidal 
$\bP_J$-equivariant local system on $C$. Let $\fS=C\cz^0_{\bP_J}$, a locally closed subvariety 
of $\bP_J$; let $cl(\fS)$ be the closure of $\fS$ in $\bP_J$. Let $\cs$ be the inverse image of
$\cs_0$ under $\fS@>>>C$ (taking unipotent part). Let $\cw=N_WW_J/W_J$ where $N_WW_J$ is the 
normalizer of $W_J$ in $W$. For any $i\in I-J$ let $N_{W_{J\cup i}}W_J$ be the normalizer of $W_J$
in $W_{J\cup i}$. From \cite{\ICC, 9.2} we see that $N_{W_{J\cup i}}W_J/W_J$ has order $2$ and 
that $\cw$ is a Coxeter group with simple reflections 
$\{\s_i;i\in I-J\}$ where $\s_i$ is the unique nonidentity element of $N_{W_{J\cup i}}W_J/W_J$, 
viewed as an element of $\cw$. Let 
$$\XX=\{(xP_J,g)\in G/P_J\T G;x\i gx\in P_J,\p_J(x\i gx)\in\fS\}.$$
Let $\bYY$ be the set of all $g\in G$ such that for some $xP_J\in G/P_J$ we have $x\i gx\in P_J$, 
$\p_J(\x\i gx)\in cl(\fS)$; this is a closed subset of $G$. Define 
$\ti\p:\XX@>>>\bYY$ by $\ti\p(xP_J,g)=g$. 
We define a local system $\hcs$ on $\XX$ by requiring that $\hcs_{(xP_J,g)}=\cs_{\p_J(x\i gx)}$.
Note that $\hcs$ is well defined by the $\bP_J$-equivariance of $\cs_0$. 
Let $K=\ti\p_!\hcs\in\cd(\bYY)$. From \cite{\ICC, 4.5, 9.2} we see that $K$ is an intersection 
cohomology complex on $\bYY$ and that we have canonically $\End(K)=\bbq[\cw]$.
Let $g\in\bYY$ and let
$$\XX_g=\{xP_J\in G/P_J;x\i gx\in P_J,\p_J(x\i gx)\in\fS\}.$$
Now $\XX_g$ is a subvariety of $\XX$ via $xP_J\m(xP_J,g)$. We denote the restriction of $\hcs$ from $\XX$ to 
$\XX_g$ again by $\hcs$. Since $H^j_c(\XX_g,\hcs)$ (for $j\in\ZZ$) is a stalk of a cohomology sheaf of $K$ at 
$g$, we see that $H^j_c(\XX_g,\hcs)$ is naturally a $\cw$-module.

We now fix a subset $J'$ of $I$ such that $J\sub J'$. We define $P_{J'}\in\cp_{J'}$ by
$P_J\sub P_{J'}$. We write $U_{J'},\p_{J'}$ instead of $U_{P_{J'}},\p_{P_{J'}}$. We shall give an alternative 
description of the restriction of the $\cw$-module $H^j_c(\XX_g,\hcs)$ to the 
subgroup $\cw'$ of $\cw$ generated by $\{\s_i;i\in J'-J\}$.

We have a commutative diagram of algebraic varieties with cartesian squares 
$$\CD \tM_{J',g}@>\tj>>\tM_{J'}@<\ti\o<<E\T X@>\tq>>X\\
@V\p''VV     @V\p'VV      @V1\T\p VV      @V\p VV\\
M_{J',g}@>j>>M_{J'}@<\o<<E\T\bY@>q>>\bY\endCD$$
where 
$$X=\{(xP_J,zU_{J'})\in(P_{J'}/P_J)\T(P_{J'}/U_{J'});x\i zx\in P_\JJ,\p_J(x\i zx)\in\fS\};$$
$\bY$ is the set of all $zU_{J'}\in P_{J'}/U_{J'}$ such that for some $xP_J\in P_{J'}/P_J$, we have 
$x\i zx\in P_J$, $\p_J(x\i zx)\in cl(\fS)$;

$\tM_{J'}$ is the set of all pairs $(xP_J,y(xU_{J'}x\i))$ where $xP_J\in G/P_J$, 
$y(xU_{J'}x\i)\in(xP_Jx\i)/(xU_{J'}x\i)$ are such that $\p_J(x\i yx)\in\fS$;

$M_{J'}$ is the set of all pairs $(xP_{J'},y(xU_{J'}x\i))$ where $xP_{J'}\in G/P_{J'}$, 
$y(xU_{J'}x\i)\in(xP_{J'}x\i)/(xU_{J'}x\i)$ are such that for some $v\in P_{J'}$, we have
$y\in xvP_Jv\i x\i$ and $\p_J(v\i x\i yxv)\in cl(\fS)$;
$$\tM_{J',g}=\{(xP_J,y(xU_{J'}x\i))\in\tM_{J'};g\in xP_Jx\i,g\i y\in xU_{J'}x\i\};$$
$$M_{J',g}=\{(xP_{J'},y(xU_{J'}x\i))\in M_{J'};g\in xP_{J'}x\i, g\i y\in xU_{J'}x\i\};$$
$$E=G/U_{J'};$$
$j,\tj$ are the obvious imbeddings; $q,\tq$ are the obvious projections;
$$\o(hU_{J'},zU_{J'})=(hP_{J'},(hzh\i)(hU_{J'}h\i)),$$ 
$$\ti\o(hU_{J'},xP_J,zU_{J'})=(hxP_J,(hzh\i)(hxU_{J'}x\i h\i));$$
$$\p(xP_J,zU_{J'})=zU_{J'},\p'(xP_J,y(xU_{J'}x\i))=(xP_{J'},y(xU_{J'}x\i)),$$ 
$$\p''(xP_J,y(xU_{J'}x\i))=(xP_{J'},y(xU_{J'}x\i)).$$
All maps in the diagram are compatible with the natural actions of $\bP_{J'}$ where the action of $\bP_{J'}$ on 
the four spaces on the left is trivial. Moreover, $\ti\o$ and $\o$ are principal $\bP_{J'}$-bundles. We define a 
local system $\bcs$ on $X$ by requiring that $\bcs_{(xP_J,zU_{J'})}=\cs_{\p_J(x\i zx)}$. We define a local system
$\dcs$ on $\tM_{J'}$ by requiring that $\dcs_{(xP_J,y(xU_{J'}x\i))}=\cs_{\p_J(x\i yx)}$. Note that $\bcs,\dcs$ 
are well defined by the $\bP_J$-equivariance of $\cs_0$. Note also that $K':=\p_!\bcs$ is like $K$ above (with 
$G$ replaced by $\bP_{J'}$) hence from \cite{\ICC, 4.5, 9.2} we see that $K'$ is an intersection cohomology sheaf on $\bY$ 
and we have canonically $\End_{\cd(\bY)}(K')=\End_{\cd_{\bP_{J'}}(\bY)}(K')=\bbq[\cw']$.

We have an isomorphism $\XX_g@>\si>>\tM_{J',g}$ given by $xP_J\m(xP_J,g(xU_{J'}x\i))$, under which these two 
varieties are identified; then the local system $\tj^*\dcs$ on $\tM_{J',g}$ becomes $\hcs$. We have 
$\tq^*\bcs=\ti\o^*\dcs$ hence $q^*\p_!\bcs=\o^*(\p'_!\dcs)$. The functors 
$$\cd_{\bP_{J'}}(\bY)@>q^*>>\cd_{\bP_{J'}}(E\T\bY)@<\o^*<<\cd(M_{J'})@>j^*>>\cd(M_{J',g})$$
induce algebra homomorphisms
$$\align&\End_{\cd_{\bP_{J'}}(\bY)}(\p_!\bcs)@>>>\End_{\cd_{\bP_{J'}}(E\T\bY)}(q^*\p_!\bcs)@<<<
\End_{\cd(M_{J'})}(\p'_!\dcs)\\&@>>>\End_{\cd(M_{J',g})}(\p''_!\tj^*\dcs)\endalign$$
of which the second one is an isomorphism since $\o$ is a principal $\bP_{J'}$-bundle. Taking the 
composition of the first homomorphism with the inverse of the second one and with the third one and identifying 
$$\End_{\cd_{\bP_{J'}}(\bY)}(\p_!\bcs)=\End_{\cd(\bY)}(\p_!\bcs)=\bbq[\cw']$$
we obtain an algebra homomorphism 
$$\bbq[\cw']@>>>\End_{\cd(M_{J',g})}(\p''_!\tj^*\dcs).$$
It follows that 
$$H^j_c(\XX_g,\hcs)=H^j_c(\tM_{J',g},\tj^*\dcs)=H^j_c(M_{J',g},\p''_!\tj^*\dcs)$$
is naturally a module over the algebra $\End_{\cd(M_{J',g})}(\p''_!\tj^*\dcs)$ hence a module over $\bbq[\cw']$.
From the definitions we see that this $\cw'$-module structure on $H^j_c(\XX_g,\hcs)$ coincides with restriction
to $\cw'$ of the $\cw$-module structure on $H^j_c(\XX_g,\hcs)$ considered above.

\subhead 2.2\endsubhead
We now assume that $i\in I-J$ and that $J'=J\cup\{i\}$. Let 
$$t_i(g)=\sum_j(-1)^j\tr(\s_i,H^j_c(\XX_g,\hcs))$$ 
where $\s_i$ acts by the $\cw$-action. By 2.1, we have
$$t_i(g)=\sum_j(-1)^j\tr(\s_i,H^j_c(M_{J',g},\p''_!\tj^*\dcs))$$
where $\s_i$ acts by the $\cw'$-action. The $\s_i$ action on $\p_!\bcs$ induce an $\s_i$-action on 
$\ch^{j'}(\p_!\bcs)$ and this induces an $\s_i$-action on $\ch^{j'}(\p_!\bcs)$ and an $\s_i$-action on 
$$H^j_c(M_{J',g},\ch^{j'}(\p''_!\tj^*\dcs)).$$ (Here $\ch^{j'}()$ denotes the $j'$-th cohomology sheaf). We have 
a spectral sequence 
$$H^j_c(M_{J',g},\ch^{j'}(\p''_!\tj^*\dcs))\imp H^{j+j'}_c(M_{J',g},\p''_!\tj^*\dcs)$$
which is compatible with the $\s_i$-actions. It follows that
$$t_i(g)=\sum_{j,j'}(-1)^{j+j'}\tr(\s_i,H^j_c(M_{J',g},\ch^{j'}(\p''_!\tj^*\dcs))).\tag a$$
We now make the further assumption that $g$ is unipotent and that $s_i$ commutes with $W_J$. 
In this case we have
$$\XX_g=\{xP_J\in G/P_J;x\i gx\in P_J,\p_J(x\i gx)\in C\};$$
moreover, the isomorphism $W_{\{i\}}\T W_J@>\si>> W_{J'}$ induced by multiplication
corresponds to a direct product decomposition $(\bP_{J'})_{ad}=H'\T H''$ where $H'\cong PGL_2(\kk)$, 
$H''\cong(\bP_J)_{ad}$.

Let $\bC$ be the image of $C$ under $\bP_J@>>>(\bP_J)_{ad}=H''$. Let $C^1$ be the unipotent class in $\bP_{J'}$ whose image in 
$(\bP_{J'})_{ad}=H'\T H''$ is $\{1\}\T\bC$. Let $C^r$ be the unipotent class in $\bP_{J'}$ whose image in 
$(\bP_{J'})_{ad}=H'\T H''$ is $(\text{regular unipotent class in $H'$})\T\bC$.
Then $C^1\cup C^r$ is exactly the set of unipotent elements in the image of $\p:X@>>>\bY$.
Now $\p:X@>>>\bY$ restricts to $\p\i(C^1)@>>>C^1$ which is a $\PP^1$-bundle and to 
$\p\i(C^r)@>>>C^r$ which is an isomorphism.
There is a well defined local system $\cs_1$ on $C^1$ whose inverse image under
$\p\i(C^1)@>>>C^1$ is the restriction of $\bcs$ to $\p\i(C^1)$.
There is a well defined local system $\cs_r$ on $C^r$ whose inverse image under
$\p\i(C^r)@>>>C^r$ is the restriction of $\bcs$ to $\p\i(C^r)$. Let
$$M_{J',g,1}=\{(xP_{J'},y(xU_{J'}x\i))\in M_{J',g};\p_{J'}(x\i gx)\in C^1\},$$
$$M_{J',g,r}=\{(xP_{J'},y(xU_{J'}x\i))\in M_{J',g};\p_{J'}(x\i gx)\in C^r\}.$$
Using the cartesian squares in 2.1 we see 
that $\p'':\tM_{J',g}@>>>M_{J',g}$ restricts to $\p''_1:\p''{}\i(M_{J',g,1})@>>>M_{J',g,1}$ which is a $\PP^1$-bundle and to
$\p''_r:\p''{}\i(M_{J',g,r})@>>>M_{J',g,r}$ which is an isomorphism. Moreover there is a well defined local system
$\cs_1$ on $M_{J',g,1}$ such that $\p''_1{}^*\cs_1$ is the restriction $\dcs_1$ of $\tj^*\dcs$ to
$\p''{}\i(M_{J',g,1})$; there is a well defined local system $\cs_r$ on $M_{J',g,r}$ such that $\p''_r{}^*\cs_r$ 
is the restriction $\dcs_r$ of $\tj^*\dcs$ to $\p''{}\i(M_{J',g,r})$.
Since $M_{J',g}=M_{J',g,1}\cup M_{J',g,r}$ is a partition and
$M_{J',g,1}$ (resp. $M_{J',g,r}$) is closed (resp. open) in $M_{J',g}$ we see that (a) implies
$$\align&t_i(g)=\sum_{j,j'}(-1)^{j+j'}\tr(\s_i,H^j_c(M_{J',g,1},\ch^{j'}(\p''_{1!}\dcs_1)))\\&+
\sum_{j,j'}(-1)^{j+j'}\tr(\s_i,H^j_c(M_{J',g,r},\ch^{j'}(\p''_{r!}\dcs_r))).\endalign$$
In the first sum over $j,j'$ we can assume that $j'\in\{0,2\}$; the $\s_i$-action is multiplication by $-1$ if $j'=2$ and by $1$ 
if $j'=0$; in the second sum over $j,j'$ we can assume that $j'=0$; moreover the $\s_i$-action is multiplication by $1$. 
(Here we use the definition of the $\s_i$ action on $K'$.) Thus, we have
$$\align&
t_i(g)=-\sum_j(-1)^j\dim H^j_c(M_{J',g,1},\ch^2(\p''_{1!}\dcs_1)))\\&
+\sum_j(-1)^j\dim H^j_c(M_{J',g,1},\ch^0(\p''_{1!}\dcs_1)))+
\sum_j(-1)^j\dim H^j_c(M_{J',g,r},\ch^0(\p''_{r!}\dcs_r)))\\&=
-\sum_j(-1)^j\dim H^j_c(M_{J',g,1},\cs_1(-2))+\sum_j(-1)^j\dim H^j_c(M_{J',g,1},\cs_1)\\&
+\sum_j(-1)^j\dim H^j_c(M_{J',g,r},\cs_r).\endalign$$
Note that the Tate twist does not affect the dimension; hence after cancellation we obtain
$$t_i(g)=\sum_j(-1)^j\dim H^j_c(M_{J',g,r},\cs_r)=
\sum_j(-1)^j\dim H^j_c(\p''{}\i(M_{J',g,r}),\dcs_r)$$
that is,
$$t_i(g)=\c(\p''{}\i(M_{J',g,r}),\dcs_r),\tag b$$
where, for an algebraic variety $X$ and a local system $\ce$ on $X$ we set
$\c(X,\ce)=\sum_j(-1)^j\dim H^j_c(X,\ce)$. (We also set $\c(X)=\c(X,\bbq)$.) 

\head 3. Computations in certain groups of semisimple rank $\le5$\endhead
\subhead 3.1 \endsubhead
Until the end of 3.14 we assume that $G_{ad}$ is of type $D_4$ and that $u\in G$ (see 1.1)
is such that $d_u=3$. Note that $u$ is unique up to conjugation.
We can write $I=\{\a,\b,\g,\o\}$ where the numbering is chosen so that each of 
$s_\a s_\o,s_\b s_\o,s_\g s_\o$ has order $3$.

The Springer representation $\r_u$ is a sum of two irreducible representations of $W$: one eight
dimensional and one six dimensional. Using 1.4(b) we see that there is a unique $S\in\un\cb_u$ such that
$J_S=\{\o\}$, a unique $\hS\in\un\cb_u$ such that $J_{\hS}=\{\a,\b,\g\}$, exactly two irreducible components
$S_{\b\g'},S_{\b'\g}$ of $\cb_u$ such that $J_{S_{\b\g'}}=J_{S_{\b'\g}}=\{\a,\o\}$, exactly two irreducible 
components $S_{\a\g},S_{\a'\g'}$ of $\cb_u$ such that $J_{S_{\a'\g'}}=J_{S_{\a'\g'}}=\{\b,\o\}$ and exactly two 
irreducible components $S_{\a\b},S_{\a'\b'}$ of $\cb_u$ such that $J_{S_{\a\b}}=J_{S_{\a'\b'}}=\{\g,\o\}$.

For any $X\in\un\cb_u$ such that $\o\in J_X$ and any $i\in\{\a,\b,\g\}$ we set
$$X_i^*=\{B\in X;\text{ any $B'$ on the same $\o$-line as $B$ is in }X_i\}.$$ 

We will show:

(a) {\it there are well defined parabolic subgroups $P^\a\ne\tP^\a$ in $\cp^\a$, $P^\b\ne\tP^\b$ in $\cp^\b$, 
$P^\g\ne\tP^\g$ in $\cp^\g$ and $P^\o$ in $\cp^\o$ such that}
$$P^\a\spa P^\o,\tP^\a\spa P^\o,P^\b\spa P^\o,\tP^\b\spa P^\o,P^\g\spa P^\o,\tP^\g\spa P^\o,$$
$$\hS=\{B\in\cb;B(\o)=P^\o\},$$
$$S_{\b\g'}=\{B\in\cb;B(\b)=P^\b,B(\g)=\tP^\g\},$$
$$S_{\b'\g}=\{B\in\cb;B(\b)=\tP^\b,B(\g)=P^\g\},$$
$$S_{\a\g}=\{B\in\cb;B(\a)=P^\a,B(\g)=P^\g\},$$
$$S_{\a'\g'}=\{B\in\cb;B(\a)=\tP^\a,B(\g)=\tP^\g\},$$
$$S_{\a\b}=\{B\in\cb;B(\a)=P^\a,B(\b)=P^\b\},$$
$$S_{\a'\b'}=\{B\in\cb;B(\a)=\tP^\a,B(\b)=\tP^\b\},$$
$$\align&S_\a^*=\{B\in\cb;B(\b)=P^\b,B(\g)=\tP^\g,B(\a)\spa P^\o\}\\&
\sqc\{B\in\cb;B(\b)=\tP^\b,B(\g)=P^\g,B(\a)\spa P^\o\},\endalign$$
$$\align&S_\b^*=\{B\in\cb;B(\a)=P^\a,B(\g)=P^\g,B(\b)\spa P^\o\}\\&
\sqc\{B\in\cb;B(\a)=\tP^\a,B(\g)\tP^\g,B(\b)\spa P^\o\},\endalign$$
$$\align&S_\g^*=\{B\in\cb;B(\a)=P^\a,B(\b)=P^\b,B(\g)\spa P^\o\}\\&
\sqc\{B\in\cb;B(\a)=\tP^\a,B(\b)=\tP^\b,B(\g)\spa P^\o\}.\endalign$$

(b) {\it Let $Y=\{B\in S;B(\o)=P^\o\}$. Then $Y\sub S_\a\cap S_\b\cap S_\g$. Moreover, $Y$ meets any $\o$-line in 
$S$ in exactly one point.}

From (b) we deduce:

(c) {\it For $i\in\{a,\b,\g\}$ we have $S_i=S_i^*\cup Y$.}
\nl
The inclusion $S_i^*\cup Y\sub S_i$ is clear.  Conversely, let $B\in S_i$ be such that $B\n Y$.
Let $L_B$ be the $i$-line containing $B$; we have $L_B\sub\cb_u$.
Let $L$ be the $\o$-line containing $B$. We have $L\sub S$ hence $L\sub\cb_u$. By (b) there exists $B'\in L$ 
such that $B'\in Y$; in particular we have $B'\in S_i$. Let $L_{B'}$ be the $i$-line containing $B'$; we have 
$L_{B'}\sub\cb_u$. We define  $\p\in\cp^{i,\o}$ by the condition that $B\sub\p$. Now $\cb^\p$ is naturally
imbedded in $\cb$; it is a flag manifold of type $A_2$ since $(s_\o s_i)^3=1$. Note that $\bar\p$ is a connected
reductive group of type $A_2$ and $\cb^\p$ can be viewed as the flag manifold of $\bar\p$. Moreover we have
$L\sub\cb^\p,L_B\sub\cb^\p,L_{B'}\sub\cb^\p$. Since $u\in B$ we have $u\in\p$ hence $Ad(u):\cb^\p@>>>\cb^\p$ is
well defined and its fixed point set contains $L\cup L_B\cup L_{B'}$. But a unipotent element in $PGL_3(\kk)$
whose fixed point set of the flag manifold contains three distinct lines must be the identity element.
In particular, for any $B''\in L$, $\Ad(u)$ acts as identity on the $i$-line through $B''$; thus the $i$-line
through $B''$ is contained in $\cb_u$. In other words we have $B\in S_i^*$. Thus we have
$S_i\sub S_i^*\cup Y$. This proves (c) (assuming (b)).

We now deduce from (a),(b) the following statement.
$$Y=S_\a\cap S_\b\cap S_\g.\tag d$$
Assume that $B\in S_\a\cap S_\b\cap S_\g$ and $B\n Y$. 
Using (c) we deduce $B\in(S^*_\a\cup Y)\cap(S^*_\b\cup Y)\cap(S^*_\g\cup Y)$. Since $B\n Y$ we deduce that
$B\in S^*_\a\cap S^*_\b\cap S^*_\g$. But from (a) we see that
$$S_\a^*\cap S_\b^*\cap S_\g^*=\emp.\tag e$$
This contradiction shows that $S_\a\cap S_\b\cap S_\g\sub Y$. The opposite inclusion is known from (b). This 
proves (d).

We now show, assuming (a):

(e) {\it $\hS\cap S_{\b\g'}$ is a single $\a$-line; $\hS\cap S_{\b'\g}$ is a single $\a$-line.}
\nl
From (a) we have
$$\align&\hS\cap S_{\b\g'}=\{B\in\cb;B(\o)=P^\o,B(\b)=P^\b,B(\g)=\tP^\g\}\\&=
\{B\in\cb;B\sub P^\o,B\sub P^\b,B\sub\tP^\g\}=\{B\in\cb;B\sub P^\o\cap P^\b\cap\tP^\g\}.\endalign$$
Since $P^\b\spa P^\o,\tP^\g\spa P^\o$ and $s_\b s_\g=s_\g s_b$, the intersection $P^\o\cap P^\b\cap\tP^\g$ is a 
parabolic subgroup in $\cp^\a$. This proves the first assertion of (e); the second assertion of (e) is proved in
 the same way.

\subhead 3.2\endsubhead
To prove  3.1(a),(b) for $G$ is the same as proving them for $G_{ad}$.
Hence we can assume that $G$ is the special orthogonal group associated to an $8$-dimensional $\kk$-vector space 
$V$ with a given nondegenerate quadratic form $Q:V@>>>\kk$ and associate symmetric bilinear form
$(,):V\T V@>>>\kk$. Until the end of 3.14 we shall adhere to this assumption.

Now, $A:=u-1:V@>>>V$ is nilpotent with Jordan blocks of sizes $3,3,1,1$.
More precisely, we can find a basis $\{e_i,e'_i;i\in\{0,1,2,3\}\}$ of $V$ such that

$Q(e_i)=Q(e'_i)=0$ for $i\in\{0,1,2,3\}$;

$(e_i,e_j)=(e'_i,e'_j)=0$ for all $i,j$;

$(e_i,e'_j)=1$ if $i\ne j$ are both odd or if $i=j$ are both even;

$(e_i,e'_j)=0$ otherwise
\nl
and such that

$Ae_0=0,Ae'_0=0$,

$Ae_1=e_2+xe_3, Ae_2=e_3, Ae_3=0$, $Ae'_1=-e'_2+x'e'_3, Ae'_2=-e'_3, Ae'_3=0$.
\nl
where $x,x'\in\kk$ satisfy $x+x'=1$.
A subspace $U$ of $V$ is said to be isotropic if $Q|_U=0$. 
Let ${}^!L_4=\spn(e_3,e'_3,e_0,e_2)$, ${}^!L'_4=\spn(e_3,e'_3,e'_0,e'_2)$.
$L_4^!=\spn(e_3,e'_3,e'_0,e_2)$, $L'_4{}^!=\spn(e_3,e'_3,e_0,e'_2)$.
These subspaces are isotropic; in 3.4 we will show that they are intrinsic to $u$; they do not depend
on the specific basis used to define $u$.
We identify $\cp^\a$ with the variety of all isotropic lines in $V$; 
$\cp^\o$ with the variety of all isotropic planes in $V$; 
$\cp^\b$  with the variety of all isotropic $4$-spaces in $V$ in the $G$-orbit of ${}^!L_4$ and ${}^!L'_4$; 
$\cp^\g$  with the variety of all isotropic $4$-spaces in $V$ in the $G$-orbit of $L_4^!$ and $L'_4{}^!$.
(In each case the identification attaches to an isotropic subspace its stabilizer in $G$.) 
We identify $\cb$ with the variety of all sequences
$(V_1,V_2,V_4,\tV_4)\in\cp^\a\T\cp^\o\T\cp^\b\T\cp^\g$
such that $V_1\sub V_2\sub V_4, V_2\sub\tV_4$.
(Such a sequence is identified with its stabilizer $\{g\in G;gV_1k=V_1,gV_2=V_2,gV_4=V_4,g\tV_4=\tV_4\}$.)

We consider the following subspaces of $V$:

$I_2=A^2V=\spn(e_3,e'_3)$, $I_4=AV=\spn(e_2,e_3,e'_2,e'_3)$, 

$K_4=\ker(A)=\spn(e_0,e_3,e'_0,e'_3)$, $K_6=\ker(A^2)=span(e_0,e_2,e_3,e'_0,e'_2,e'_3)$.
\nl
We have 
$I_2\sub I_4\sub K_6$, $I_2\sub K_4\sub K_6$, $K_4\cap I_4=I_2$.

\subhead 3.3\endsubhead
We show:

(a) {\it There are exactly two subspaces $V_4\in\cp^\b$ such that $I_2\sub V_4$ and $\dim(AV_4)\le1$. They are
${}^!L_4,{}^!L'_4$.}

(b) {\it There are exactly two subspaces $V_4\in\cp^\g$ such that $I_2\sub V_4$ and $\dim(AV_4)\le1$. They are
$L_4^!,L'_4{}^!$.}
\nl
We prove (a). The subspaces $V_4\in\cp^\b$ such that $I_2\sub V_4$ are determined by their intersection with
$L_4^!,L'_4{}^!$; this intersection is of the form $\spn(e_3,e'_3,ae_0+be'_2)$ where $(a,b)\in\kk^2-\{0,0\}$
hence the required $4$-subspaces are of the form $\spn(e_3,e'_3,ae_0+be'_2,xe_0+x'e'_0+ze_2+z'e'_2)$ where 
$ae_0+be'_2,xe_0+x'e'_0+ze_2+z'e'_2$ are linearly independent and we have $ax'+bx=0,xx'+zz'=0$. Moreover the
condition that $\dim(AV_4)\le1$ is that $\dim\spn(-be'_3,ze_3-z'e'_3)\le1$, that is $bz=0$.
Assume first that $a\ne,b\ne0$; then $z=0$, $xx'=0$. From $ax'+bx=0,xx'=0$ we deduce that $x=x'=0$.
Since $ae_0+be'_2,z'e'_2$ are linearly independent we see that $z'\ne0$ and our $V_4$ is
$\spn(e_3,e'_3,ae_0+be'_2,e'_2)=\spn(e_3,e'_3,e_0,e'_2)\in\cp^\g$, a contradiction.
Assume now that $a=0,b\ne0$; then $z=0,x=0$.
Since $be'_2,x'e'_0+z'e'_2$ are linearly independent we see that $x'\ne0$ and our $V_4$ is
$\spn(e_3,e'_3,e'_2,x'e'_0+z'e'_2)=\spn(e_3,e'_3,e'_2,e'_0)={}^!L'_4$.
Assume now that $a\ne0,b=0$; then $x'=0$ and $zz'=0$.
Since $ae_0,xe_0+ze_2+z'e'_2$ are linearly independent we have either $z=0,z'\ne0$ or $z\ne0,z'=0$.
If $z=0,z'\ne0$ then our $V_4$ is
$\spn(e_3,e'_3,e_0,xe_0+z'e'_2)=\spn(e_3,e'_3,e_0,e'_2)\in\cp^\g$, a contradiction.
If $z\ne0,z'=0$ then our $V_4$ is
$\spn(e_3,e'_3,e_0,xe_0+ze_2)=\spn(e_3,e'_3,e_0,e_2)={}^!L_4$.
Thus $V_4\in\{{}^!L_4,{}^!L'_4\}$. Conversely it is clear that if $V_4\in\{{}^!L_4,{}^!L'_4\}$ then $V_4$
satisfies the requirements of (a). This proves (a).

The proof of (b) is entirely similar to that of (a).

\subhead 3.4\endsubhead
We set $L_1=\spn(e_3)\in\cp^\a$, $L'_1=\spn(e'_3)\in\cp^\a$.  We have $A({}^!L_4)=A(L_4^!)=L_1$, 
$A({}^!L'_4)=A(L'_4{}^!)=L'_1$. In particular, if $V_4$ is as in 3.3(a),(b), then $AV_4\in\{L_1,L'_1\}$ is a one 
dimensional (isotropic) subspace of $I_2$.

From 3.3(a),(b) we see that ${}^!L_4,L_4^!,{}^!L'_4,L'_4{}^!$ are intrinsic to $u$ and do not depend on the
specific basis used to define $u$. It follows also that $L_1,L'_1$ are intrinsic to $u$.
Note that $L_1,L'_1$ are the two lines in $I_2$ which are isotropic for the quadratic form
$I_2@>>>\kk$ given by $v\m Q(\ti v)$ where $\ti v$ is any vector in $I_4$ such that $A\ti v=v$. 

\subhead 3.5\endsubhead
Note that $Q$ induces a nondegenerate quadratic form on $K_4/I_2$ which has exactly two isotropic lines. Hence 
there are exactly two isotropic $3$-spaces contained in $K_4$ and containing $I_2$. They are 
$\spn(e_3,e'_3,e_0)={}^!L_4\cap L'_4{}^!$ and $\spn(e_3,e'_3,e'_0)={}^!L'_4\cap L_4^!$.

\subhead 3.6\endsubhead
Let 
$$\tY=\{(V_1,V_2,V_4,\tV_4)\in\cb;V_2=I_2\}.$$
Let $(V_1,V_2,V_4,\tV_4)\in\tY$. We have clearly $AV_1=0$. Moreover, since $V_4$ is an isotropic subspace
containing $I_2$ we must have $V_4\sub K_6=\ker(A^2)$ hence $AV_4\sub \ker(A)\cap AV=I_4\cap K_4=I_2\sub V_4$. 
Thus, $AV_4\sub V_4$. Similarly we have $A\tV_4\sub\tV_4$. We see that $(V_1,V_2,V_4,\tV_4)\in\cb_u$. Thus,
$\tY\sub\cb_u$. Note that $\tY$ is isomorphic to $\PP^1\T\PP^1\T\PP^1$ hence it is a closed irreducible 
subvariety of $\cb_u$ of dimension $3$. Thus $\tY$ is an irreducible component of $\cb_u$. It satisfies 
$J_{\tY}=\{\a,\b,\g\}$. Hence $\tY=\hS$.

\subhead 3.7\endsubhead
Let
$$N_1=\{(V_1,V_2,V_4,\tV_4)\in\cb;V_4={}^!L_4,\tV_4=L'_4{}^!\},$$
$$N_2=\{(V_1,V_2,V_4,\tV_4)\in\cb;V_4={}^!L'_4,\tV_4=L_4^!\}.$$
We have $A({}^!L_4\cap L'_4{}^!)=0$, $A({}^!L'_4\cap L_4^!)=0$ hence if 
$(V_1,V_2,V_4,\tV_4)$ is in $N_1$ or $N_2$
then $AV_2=0$. Moreover, each of ${}^!L_4,L'_4{}^!,{}^!L'_4,L_4^!$ is $A$-stable. We see that 
$(V_1,V_2,V_4,\tV_4)\in\cb_u$. Thus, $N_1\sub\cb_u$, $N_2\sub\cb_u$. Note that $N_1$ are $N_2$ are isomorphic to 
the flag manifold of $GL_3(\kk)$ hence they are closed irreducible subvarieties of $\cb_u$ of dimension $3$. Thus
they are (distinct) irreducible components of $\cb_u$. They satisfy $J_{N_1}=J_{N_2}=\{\a,\o\}$.
Hence $N_1,N_2$ are the same as $S_{\b\g'},S_{\b'\g}$ (up to order).

\subhead 3.8\endsubhead
Let
$$N_3=\{(V_1,V_2,V_4,\tV_4)\in\cb;V_1=L_1, \tV_4=L_4^!\},$$
$$N_4=\{(V_1,V_2,V_4,\tV_4)\in\cb;V_1=L'_1, \tV_4=L'_4{}^!\}.$$
Let $(V_1,V_2,V_4,\tV_4)\in N_3$. Let $V_3=V_4\cap\tV_4$. Since $A(\tV_4)=L_1$ we have 
$A(V_3)\sub L_1=V_1\sub V_3$ that is $uV_3=V_3$; moreover, $A(V_2)\sub L_1=V_1$ hence $AV_2\sub V_2$.
Now $V_4$ is the only subspace in its $G$-orbit that
contains $V_3$; since $uV_4\supset uV_3=V_3$ it follows that $uV_4=V_4$. The inclusions $AV_1\sub V_1$,
$A(\tV_4)\sub\tV_4$ are obvious. We see that $(V_1,V_2,V_4,\tV_4)\in\cb_u$. Thus, $N_3\sub\cb_u$. Note that 
$N_3$ is isomorphic to the space of pairs $(V_2,V_3)$ where $V_2\sub V_3$ are subspaces of $\tV_4$ with
$\dim V_2=2,\dim V_3=3$; hence $N_3$ is isomorphic to the the flag manifold of $GL_3(\kk)$. Thus
$N_3$ is a closed irreducible subvariety of $\cb_u$ of dimension $3$. Thus
$N_3\in\un\cb_u$. We have  $J_{N_3}=\{\b,\o\}$. Similarly, $N_4\in\un\cb_u$ and $J_{N_4}=\{\b,\o\}$. 
Hence $N_3,N_4$ are the same as $S_{\a\g},S_{\a'\g'}$ (up to order).

\subhead 3.9\endsubhead
Let
$$N_5=\{(V_1,V_2,V_4,\tV_4)\in\cb;V_1=L_1, V_4={}^!L_4\},$$
$$N_6=\{(V_1,V_2,V_4,\tV_4)\in\cb;V_1=L'_1, V_4={}^!L'_4\}.$$
A proof completely similar to that in 3.8 shows that
$N_5\in\un\cb_u$, $N_6\in\un\cb_u$. We have $J_{N_5}=J_{N_6}=\{\g,\o\}$.  
Hence $N_5,N_6$ are the same as $S_{\a\b},S_{\a'\b'}$ (up to order).

\subhead 3.10\endsubhead
Let $X=\{(V_1,V_2,V_4,\tV_4)\in\cb;V_1\sub I_2\sub V_4\cap\tV_4,A(V_4\cap\tV_4)\sub V_1\}$. We show:

(a) {\it $X$ is a smooth irreducible projective variety of dimension $3$.}
\nl
The fact that $X$ is a  projective variety is obvious. Now let $Z$ be the variety of all pairs $V_1,V_3$ of 
isotropic subspaces of $V$ of dimension $1$ and $3$ respectively such that 
$V_1\sub I_2\sub V_3$ and $AV_3\sub V_1$. Clearly $\z:(V_1,V_2,V_4,\tV_4)\m(V_1,V_4\cap\tV_4)$ makes $X$
into a $\PP^1$-bundle over $Z$. Hence it is enough to show that $Z$ is a smooth irreducible surface.
Now the space of lines in $I_2$ is a projective line and the space of $V_3$ containing $I_2$ and contained in
$K_6$ (without the isotropy condition) is a projective $3$-space. Hence we can identify $Z$ with
$$\{(z_0,z'_0,z_2,z'_2),(y,y'))\in\PP^3\T \PP^1;z_0z'_0+z_2z'_2=0, z_2y'+z'_2y=0\}.$$

The subset of $Z$ where $z'_2\ne0, y'\ne0$ is
$$\align&\{(z_0,z'_0,z_2;y)\in \kk^3\T\kk;z_0z'_0+z_2=0, z_2+y=0\}\\&=
\{(z_0,z'_0,z_2)\in\kk^3;z_0z'_0+z_2=0\}\cong\kk^2.\endalign$$
Similarly the subset of $Z$ where $z_2\ne0, y\ne0$ is $\cong\kk^2$. The subset where $z_0\ne0, y\ne0$ is
$$\{(z'_0,z_2,z'_2,y')\in\kk^4;z'_0+z_2z'_2=0, z_2y'+z'_2=0\}\cong\kk^2.$$
The subset where $z_0\ne0, y'\ne0$ is $\cong\kk^2$; the subset where $z'_0\ne0, y\ne0$ is $\cong\kk^2$;
the subset where $z'_0\ne0, y'\ne0$ is $\cong\kk^2$.
Thus, $Z$ is a union of $8$ smooth irreducible surfaces (open in $Z$) hence $Z$ is a smooth surface. 
Since the $8$ open subsets above have a nonempty intersection, we see that $Z$ is irreducible.
This proves (a).

Now if $(V_1,V_2,V_4,\tV_4)\in X$ then setting $V_3=V_4\cap\tV_4$ we have
$AV_1=0$ (since $V_1\sub I_2$), $A(V_3)\sub V_1\sub V_3$, $AV_2\sub V_1\sub V_2$ hence  $uV_1\sub V_1$, 
$uV_3\sub V_3$, $uV_2\sub V_2$; since $V_4\in\cp^\b,\tV_4\in\cp^\g$ are determined uniquely
by the conditions $V_3\sub V_4$, $V_3\sub\tV_4$, we have necessarily $uV_4=V_4$, $u\tV_4=\tV_4$. Thus,
$X\sub\cb_u$. Using this and (a) together with $\dim\cb_u=3$ we see that 

(b) {\it $X$ is an irreducible component of $\cb_u$.}
\nl
It follows that the subsets $X_\a,X_\b,X_\g$ of $X$ are well defined (see 1.6). 

We show:

(c) $\c(X)=12$.
\nl
Here $\c()$ is as in 2.2. Since $X$ is a $\PP^1$-bundle over $Z$, an equivalent statement is:

(d) $\c(Z)=6$.
\nl
The space $\fU$ of isotropic $3$-spaces in $V$ that contain $I_2$ hence are contained in $K_6$ can be identified 
with the space of isotropic lines in $K_6/I_2$ with its obvious quadratic form hence it is a product 
$\PP^1\T\PP^1$. 
The map $Z@>>>\fU$, $(V_1,V_3)\m V_3$ is an isomorphism over the complement of two points in $\fU$, namely 
${}^!L_4\cap L'_4{}^!$ and ${}^!L'_4\cap L_4^!$ and has fibres $\PP^1$  at each of those two points. It follows 
that $\c(Z)=\c(\PP^1\T \PP^1)-2+2\c(\PP^1)=4-2+4=6$, as desired.

\subhead 3.11\endsubhead
Clearly, $X$ is a union of $\o$-lines. Hence for $i\in\{\a,\b,\g\}$, $X_i^*$ is defined as in 3.1. We have 
$$X_\a=\{(V_1,V_2,V_4,\tV_4)\in X;\text{ for any line $V'_1$ in $V_2$ we have }AV'_1=0\}$$
hence
$$X_\a=\{(V_1,V_2,V_4,\tV_4)\in X;V_2\sub K_4\}.$$
We deduce that
$$X_\a^*=\{(V_1,V_2,V_4,\tV_4)\in X;\text{ for any $V'_2$ with $V_1\sub V'_2\sub V_4\cap\tV_4$ we have }
V'_2\sub K_4\}.$$
Since the various $V'_2$ such that $V_1\sub V'_2\sub V_4\cap\tV_4$ generate $V_4\cap\tV_4$ we see that
$$X_\a^*=\{(V_1,V_2,V_4,\tV_4)\in X; I_2\sub V_4\cap\tV_4\sub K_4\}$$
If $(V_1,V_2,V_4,\tV_4)\in X_\a^*$ then $A(V_4\cap\tV_4)=0$ hence
$\dim A(V_4)\le 1$, $\dim A(\tV_4)\le 1$ (since $V_4\cap\tV_4$ has codimension $1$ in $V_4$ and
in $\tV_4$). Using 3.3(a),(b) we deduce that $V_4\in\{{}^!L_4,{}^!L'_4\}$,
$\tV_4\in\{L_4^!,L'_4{}^!\}$. Thus we have either
$V_4={}^!L_4,\tV_4=L'_4{}^!$ or $V_4={}^!L'_4,\tV_4=L_4^!$. We see that
$$\align&X_\a^*=\{(V_1,V_2,V_4,\tV_4)\in X;V_4={}^!L_4,\tV_4=L'_4{}^!\}\\&\sqc
\{(V_1,V_2,V_4,\tV_4)\in X;V_4={}^!L'_4,\tV_4=L_4^!\}.\endalign$$
(The right hand side is clearly contained in the left hand side.)
Since $A({}^!L_4\cap L'_4{}^!)=0,A({}^!L'_4\cap L_4^!)=0$, we have
$$\align&X_\a^*=\{(V_1,V_2,V_4,\tV_4)\in\cb; V_1\sub I_2,V_4={}^!L_4,\tV_4=L'_4{}^!\}\\&\sqc
\{(V_1,V_2,V_4,\tV_4)\in\cb; V_1\sub I_2,V_4={}^!L'_4,\tV_4=L_4^!\}.\tag a\endalign$$

\subhead 3.12\endsubhead
We have
$$X_\b=\{(V_1,V_2,V_4,\tV_4)\in X;\text{ for any  $V'_4\in\cp^\b$ with $V_2\sub V'_4$ we have }u(V'_4)=V'_4\}.$$
Let $(V_1,V_2,V_4,\tV_4)\in X_\b$. Then for any $V'_4\in\cp^\b$ with $V_2\sub V'_4$ we have that
$\tV_4\cap V'_4$ is an isotropic $3$-space containing $V_2$; it is $u$-stable (since $\tV_4$ and $V'_4$ are
$u$-stable). Moreover all isotropic $3$-spaces containing $V_2$ and contained in
$\tV_4$ are obtained in this way hence are all $u$-stable; it follows that $u$ acts as $1$ on $\tV_4/V_2$
that is, $A(tV_4)\sub V_2$. Thus, $X_\b\sub\{(V_1,V_2,V_4,\tV_4)\in X;A(\tV_4)\sub V_2\}.$ A similar argument
shows the reverse inclusion. Thus
$$X_\b=\{(V_1,V_2,V_4,\tV_4)\in X;A(\tV_4)\sub V_2\}.$$
We deduce that
$$\align&X_\b^*=\{(V_1,V_2,V_4,\tV_4)\in X;\text{ for any $V'_2$ with 
$V_1\sub V'_2\sub V_4\cap\tV_4$ we have }\\& A(\tV_4)\sub V'_2\}.\endalign$$
Since the intersection of all $V'_2$ such that $V_1\sub V'_2\sub V_4\cap\tV_4$ is $V_1$ we see that
$$X_\b^*=\{(V_1,V_2,V_4,\tV_4)\in X;A(\tV_4)\sub V_1\}$$
that is,
$$X_\b^*=\{(V_1,V_2,V_4,\tV_4)\in\cb;V_1\sub I_2\sub V_4\cap\tV_4, A(\tV_4)\sub V_1\}.$$
Using 3.3(b), we see that for $\tV_4$ in the right hand side we have $\tV_4\in\{L_4^!,L'_4{}^!\}$ so that $V_1$ is 
necessarily $L_1$ (if $\tV_4=L_4^!$) or $L'_1$ (if $\tV_4=L'_4{}^!$). Thus,
$$\align&X_\b^*=\{(V_1,V_2,V_4,\tV_4)\in\cb;I_2\sub V_4,V_1=L_1,\tV_4=L_4^!\}\\&\sqc
\{(V_1,V_2,V_4,\tV_4)\in\cb;I_2\sub V_4,V_1=L'_1,\tV_4=L'_4{}^!\}.\tag a\endalign$$
An entirely similar argument yields:
$$X_\g=\{(V_1,V_2,V_4,\tV_4)\in S;A(V_4)\sub V_2\},$$
$$\align&X_\g^*=\{(V_1,V_2,V_4,\tV_4)\in\cb;I_2\sub\tV_4,V_1=L_1,V_4={}^!L_4\}\\&
\sqc\{(V_1,V_2,V_4,\tV_4)\in\cb;I_2\sub\tV_4,V_1=L'_1,V_4={}^!L'_4\}.\tag b\endalign$$
From 3.11(a) we see that 
$X_\a^*$ has two irreducible components, each one being a $\PP^1$-bundle over $\PP^1$ hence is
two-dimensional; we see that $X_\a^*\ne X$ hence $X$ is not a union of $\a$-lines.
Similarly from (a),(b)
 we see that $X_\b^*\ne X$, $X_\g^*\ne X$ hence $X$ is not a union of $\b$-lines and
$X$ is not a union of $\g$-lines. Thus we have $J_X=\{\o\}$. It follows that 
$$X=S.\tag c$$

\subhead 3.13\endsubhead
We have
$$\align&Y=\{(V_1,V_2,V_4,\tV_4)\in X;V_2=I_2\}\\&
=\{(V_1,V_2,V_4,\tV_4)\in\cb;V_2=I_2,A(V_4\cap\tV_4)\sub V_1\}.\endalign$$
Assume that $(V_1,V_2,V_4,\tV_4)\in Y$. If $V'_1$ is an isotropic line contained in $I_2$ then $AV'_1=0$ (since 
$AI_2=0$). This shows that $Y\sub X_\a=S_\a$. If $U$ is an isotropic $4$-space containing $I_2$ then 
$U\sub K_6=\ker(A^2)$ hence $AU\sub\ker(A)\cap AV=K_4\cap I_4=I_2\sub U$; thus, $AU\sub U$. This shows that 
$Y\sub X_\b=S_\b$ and
$Y\sub X_\g=S_\g$. This proves the first assertion of 3.1(b). 
Now let $(V_1,V_2,V_4,\tV_4)\in X$ and let
$L$ be the $\o$-line in $X$ containing $(V_1,V_2,V_4,\tV_4)$. By the definition of $X$ if $V_2$ is replaced by 
$I_2$ the resulting quadruple $(V_1,I_2,V_4,\tV_4)$ belongs to $X$ (and even to $Y$). This proves the second 
assertion of 3.1(b).

\subhead 3.14\endsubhead
From 3.11(a) and 3.12(a),(b),(c), we see that 3.1(a) holds. Here $P^\a$ (resp. $\tP^\a$) is 
the stabilizer of $L_1$ (resp. $L'_1$); $P^\b$ (resp. $\tP^\b$) is the stabilizer of ${}^!L_4$ 
(resp. ${}^!L'_4$); $P^\g$ (resp. $\tP^\g$) is the stabilizer of $L_4^!$ (resp. $L'_4{}^!$).  

\subhead 3.15\endsubhead
Until the end of 3.18 we assume that $G_{ad}$ is of type $D_5$.
In this case we can write $I=\{\a,\b,\g,\d,\o\}$ where the numbering is chosen so that each of
$s_\a s_\o,s_\b s_\o, s_\g s_\o,s_\a s_\d$ has order $3$.
We shall assume that $u\in G$ (see 1.1) 
is such that $d_u=3$; this determines $u$ uniquely up to 
conjugacy. We can find
$P\in\cp^\d$ such that $u\in P$ and $\bu:=\p_P(u)\in\bP$ satisfies $d_{\bu}=3$ and, if $C$ is the 
conjugacy class of $u$ in $G$, then $C\cap\p_P\i(\bu)$ is open dense in $\p_P\i(\bu)$.
Since the adjoint group of $\bP$ is of type $D_4$, the irreducible component $S$ of $\cb^{\bP}_{\bu}$ 
is defined as in 3.1 with $G,u$ replaced by $\bP,\bu$. 
Note that the Weyl group of $\bP$ is canonically identified with the subgroup of $W$ generated by
$\{s_\a,s_\b,s_\g,s_\o\}$. The imbedding $\cb^{\bP}@>>>\cb$, $B_1\m\p_P\i(B_1)$ restricts to an
imbedding $\cb^{\bP}_{\bu}@>>>\cb_u$ and to an isomorphism of $S$ onto an irreducible component 
$\tS$ of $\cb_u$. (See 1.5.) This imbedding carries any $i$-line in $\cb^{\bP}$ (where $i\in\{\a,\b,\g,\o\}$) to
an $i$-line in $\cb$. 

Using 1.4(b), we see that there is a unique $E\in\un\cb_u$ such that $J_E=\{\b,\g,\d\}$.

The inverse images under $\p_P$ of the parabolic subgroups 

$P^\a,\tP^\a,P^\b,\tP^\b,P^\g,\tP^\g,P^\o$ 
\nl
as in 3.1(a) (with $G,u$ 
replaced by $\bP,\bu$) are denoted again by the same letters; thus we have 

$P^\a\ne\tP^\a$ in $\cp^\a$, $P^\b\ne\tP^\b$ in $\cp^\b$, $P^\g\ne\tP^\g$ in $\cp^\g$, 
$P^\o\in\cp^\o$, 
\nl
where $\cp^i$ refers to $G$.

\subhead 3.16\endsubhead
We will show that there is a well defined parabolic subgroup $\hP^\a\in\cp^\a$ such that

(a) $\hP^\a\n\{P^\a,\tP^\a\}$, $\hP^\a\spa P^\o$ and

(b) $\tS_\d=\{B\in\tS;B(\a)=\hP^\a\}$,

(c) $E=\{B\in\cb;B(\a)=\hP^\a,B(\o)=P^\o\}$.
\nl
To prove this is the same as proving it for $G_{ad}$. 
Hence we can assume that $G$ is the special orthogonal group associated to a $10$-dimensional $\kk$-vector space 
$V$ with a given nondegenerate quadratic form $Q:V@>>>\kk$ and associate symmetric bilinear form
$(,):V\T V@>>>\kk$. Until the end of 3.18 we shall adhere to this assumption.
Now $A:=u-1:V@>>>V$ is nilpotent. If $p\ne2$, $A$ has Jordan blocks of sizes $5,3,1,1$.
If $p=2$, $A$ has Jordan blocks of sizes $4,4,1,1$ and moreover we have $Q(A^2x)\ne0$ for some $x\in V$. 
These conditions describe completely the conjugacy class of $u$. 

We identify $\cp^\d$ with the variety of all isotropic lines in $V$; $\cp^\a$ with the variety of all isotropic 
planes in $V$; $\cp^\o$ with the variety of all isotropic $3$-spaces in $V$. (In each case the identification 
attaches to an isotropic subspace its stabilizer in $G$.) We will describe later $\cp^\b$ and $\cp^\g$.

We can find a basis $\{e_i,e'_i;i\in\{0,1,2,3\}\}\sqc\{f,f'\}$ of $V$ such that

$P$ is the stabilizer in $G$ of the line $\kk f$;

$Q(e_i)=Q(e'_i)=0$ for $i\in\{0,1,2,3\}$, $Q(f)=Q(f')=0$;

$(e_i,e_j)=(e'_i,e'_j)=0$ for all $i,j$ and $(f,e_i)=(f',e_i)=0$ for all $i$;

$(e_i,e'_j)=1$ if $i\ne j$ are both odd or if $i=j$ are both even and $(f,f')=1$;

$(e_i,e'_j)=0$ otherwise,
\nl
and such that

$Af'=a_0e_0+a_1e_1+a_2e_2+a_3e_3+a'_0e'_0+a'_1e'_1+a'_2e'_2+a'_3e'_3+df$

$Ae_0=b_0f, Ae'_0=b'_0f$,

$Ae_1=e_2+xe_3+b_1f, Ae_2=e_3+b_2f, Ae_3=b_3f$, $Ae'_1=-e'_2+x'e'_3+b'_1f$, 

$Ae'_2=-e'_3+b'_2f, Ae'_3=b'_3f,Af=0$,
\nl
where $x,x'\in\kk$ satisfy $x+x'=1$, compare with 3.2. Here 

$a_0,a_1,a_2,a_3,a'_0,a'_1,a'_2,a'_3,b_0,b_1,b_2,b_3,b'_0,b'_1,b'_2,b'_3,d$ 
\nl
are elements of $\kk$ such that

$d+a_0a'_0+a_1a'_3+a_2a'_2+a_3a'_1=0$, $a'_0+b_0=0, a_0+b'_0=0$,

$a'_3+a'_2+xa'_1+b_1=0,   a_3-a_2+x'a_1+b'_1=0$,
$a'_2+a'_1+b_2=0,  a_2-a_1+b'_2=0$,
              
$a'_1+b_3=0, a_1+b'_3=0$.
\nl
(These equations express the fact that $u$ is an isometry for $Q$.)

We identify $\cp^\b$ with the variety of all isotropic $5$-spaces in $V$ in the $G$-orbit of
$\spn(f,e_3,e'_3,e_0,e_2)$; $\cp^\g$ with the variety of all isotropic $5$-spaces in $V$ in the $G$-orbit of
$\spn(f,e_3,e'_3,e'_0,e_2)$.
(In each case the identification attaches to an isotropic subspace its stabilizer in $G$.) 
We identify $\cb$ with the variety of all sequences
$(V_1,V_2,V_3,V_5,\tV_5)\in\cp^\d\T\cp^\a\T\cp^\o\T\cp^\b\T\cp^\g$
such that $V_1\sub V_2\sub V_3\sub V_5, V_3\sub\tV_5$.
(Such a sequence is identified with its stabilizer 
$\{g\in G;gV_1=V_1,gV_2=V_2,gV_3=V_3,gV_5=V_5,g\tV_5=\tV_5\}$.)

We now compute the powers $A^k$ for $k=2,3,4$ (here $*$ denotes an element of $\kk$).

$A^2f'=a_1e_2+(a_2+xa_1)e_3-a'_1e'_2+(-a'_2+x'a'_1)e'_3+ *f$,

$A^2e_0=0$, $A^2e_1=e_3+*f$, $A^2e_2=-a'_1f$, $A^2e_3=0$, 

$A^2e'_0=0$, $A^2e'_1=e'_3+*f$, $A^2e'_2=-a_1f$, $A^2e'_3= 0$, $A^2f=0$;

$A^3f'=a_1e_3+a'_1e'_3+*f$, $A^3e_0=0$, $A^3e_1=-a'_1f$, $A^3e_2=0$, $A^3e_3=0$, $A^3e'_0=0$,

$A^3e'_1=-a_1f$, $A^3e'_2=0$, $A^3e'_3=0$, $A^3f=0$;

$A^4f'=-2a_1a'_1f$, $A^4e_i=A^4e'_i=0$ for $i=0,1,2,3$, $A^4f=0$.
\nl
We see that $A^4=0$ if $p=2$ and $A^5=0$ without restriction on $p$. Since $A$ has a Jordan block of size $5$ 
(if $p\ne2$) and one of size $4$ (if $p=2$) we see that $A^4\ne0$ (if $p\ne2$) and $A^3\ne0$ (if $p=2$). It 
follows that $a_1a'_1\ne0$ if $p\ne2$ and $(a_1,a'_1)\ne(0,0)$ in any case.

Let $\ci$ be the image of the map $V@>>>\kk$, $v\m Q(A^2v)$. Since 
$$A^2V=\spn(f,e_3,e'_3,a_1e_2-a'_1e'_2),$$ 
$\ci$ is the same as the image of the map $\kk^4@>>>\kk$, 
$$(x,y,x,t)\m Q(xf+ye_3+ze'_3+t(a_1e_2-a'_1e'_2)=Q(t(a_1e_2-a'_1e'_2))=-t^2a_1a'_1.$$
Thus, if $p\ne2$ we have $\ci=\kk$ since $a_1a'_1\ne0$; if $p=2$ then we already know that $\ci\ne0$ hence we 
must have $a_1a'_1\ne0$. Thus, without assumption on $p$ we have $a_1a'_1\ne0$ and $\ci=\kk$.

\subhead 3.17\endsubhead
Let $I_3=\spn(f,e_3,e'_3)\in\cp^\o$, $I'_2=\spn(f,a_1e_3-a'_1e'_3)$,
$K_7=\{v\in V;(v,V_3)=0\}=\spn(f,e_3,e'_3,e_2,e'_2,e_0,e'_0)$.

Note that $I'_2=I_3\cap\ker(A)$; moreover $I'_2,I_3$ are isotropic subspaces.
From the definitions we have
$$\tS=\{(V_1,V_2,V_3,V_5,\tV_5)\in\cb;V_1=\kk f\sub V_2\sub I_3\sub V_5\cap\tV_5,A(V_5\cap\tV_5)\sub V_2\}.$$
From the definitions, $\tS_\d$ is the set of all $(V_1,V_2,V_3,V_5,\tV_5)\in\tS$ such that for any isotropic line
$V'_1$ contained in $V_2$ we have $AV'_1\sub V'_1$ that is, $AV'_1=0$. Since $V_2$ is generated by such $V'_1$, 
we see that
$$\align&\tS_\d=\{(V_1,V_2,V_3,V_5,\tV_5)\in\tS;AV_2=0\}\\&=
\{(V_1,V_2,V_3,V_5,\tV_5)\in\tS;V_2\sub I_3\cap \ker(A)\}
=\{(V_1,V_2,V_3,V_5,\tV_5)\in\tS;V_2=I'_2\}.\endalign$$
This shows that 3.16(b) holds where $\hP^\a$ is the parabolic subgroup in $\cp^\a$ which is the stabilizer of 
$I'_2$. From the definitions we see that $P^\a$ (resp. $\tP^\a$) is the stabilizer of $\spn(e_3,f)$
(resp. of $\spn(e'_3,f)$). Since $\spn(e_3,f),\spn(e'_3,f),\spn (a_1e_3-a'_1e_3,f)$ are distinct
(recall that $a_1a'_1\ne0$) we see that $\hP^\a\n\{P^\a,\tP^\a\}$. 
Since $P^\o$ is the stabilizer of $I_3$ and $I'_2\sub I_3$ we see that $\hP^\a\spa P^\o$.
Thus, 3.16(a) holds.

\subhead 3.18\endsubhead
Let 
$$E'=\{(V_1,V_2,V_3,V_5,\tV_5)\in\cb;V_2=I'_2,V_3=I_3\}.$$
If $(V_1,V_2,V_3,V_5,\tV_5)\in E'$ then $AV_1=AV_2=AV_3=0$; moreover we have 
$V_5\sub K_7$ hence $AV_5\sub AK_7\sub I_3\sub V_5$. Similarly, $A\tV_5\sub I_3\sub V_5$. Thus we have
$(V_1,V_2,V_3,V_5,\tV_5)\in\cb_u$. We see that $E'\sub\cb_u$. Now $E'$ is isomorphic to 
$\PP^1\T\PP^1\T\PP^1$ hence it is a closed irreducible $3$-dimensional subvariety of $\cb_u$. Thus $E'$ is
an irreducible component of $\cb_u$. We have $J_{E'}=\{\d,\b,\g\}$. We see that $E'=E$. From the description
of $E=E'$ given above we see that 3.16(c) holds.

\subhead 3.19\endsubhead
Until the end of 3.23 we assume that $G_{ad}$ is of type $A_5$ and that $u\in G$ (see 1.1)
is such that
$d_u=3$ and the Springer representation $\r_u$ is the irreducible representation of $W$ of dimension five 
appearing in the third symmetric power of the reflection representation of $W$. 
 Note that $u$ is unique up to conjugation. We can write $I=\{1,2,3,4,5\}$ where the numbering is chosen 
so that each of $s_1s_2,s_2s_3,s_3s_4,s_4s_5$ has order $3$.

Using  1.4(b) we see that there is a unique 
$T\in\un\cb_u$ such that $J_T=\{3\}$ and a unique $M\in\un\cb_u$ such that $J_M=\{1,3,5\}$. We shall prove

(a) {\it there are well defined parabolic subgroups $P^2\in\cp^2$, $P^4\in\cp^4$ such that
$M=\{B\in\cb;B(2)=P^2,B(4)=P^4\}$;

we have $T_1=T_5=\{B\in T; B(2)=P^2,B(4)=P^4\}=M\cap T$ (we denote it by $T_{15}$); this is a 
$\PP^1$-bundle over $\PP^1$ and is a union of $3$-lines;

we have $T_2=T_4$ (we denote it by $T_{24}$); it intersects any $3$-line in $T$ in exactly one point; 

the intersection $T_{15}\cap T_{24}$ is isomorphic to $\PP^1$;

there is a unique morphism $\vt:T_{24}@>>>T_{15}\cap T_{24}$ such that for any $B\in T_{24}$ we have 
$(B,\vt(B))\in\co_w$ for some $w\in W_{24}$; moreover, $\vt$ is a $\PP^1$-bundle.}
\nl
From (a) we see that 

(b) $\c(T_{24})=4$, $\c(T)=8$.
\nl
where $\c()$ is as in 2.2.

\subhead 3.20\endsubhead
To prove  3.19(a) for $G$ is the same as proving it for $G_{ad}$. Hence we can assume that 
$G=GL(V)$ where $V$ is 
a $6$-dimensional $\kk$-vector space. Until the end of 3.23 we shall adhere to this assumption.
We can write $u=A+1$ where $A:V@>>>V$ is nilpotent with Jordan blocks of sizes $3,3$. Let
$K_2=\ker A=A^2(V)$, $K_4=\ker A^2=A(V)$. We have $K_2\sub K_4$.

For $i\in I$, we identify $\cp^i$ with the variety of all $i$-dimensional subspaces of $V$.
(The identification attaches to a subspace its stabilizer in $G$.) 
We identify $\cb$ with the variety of all sequences $(V_1,V_2,V_3,V_4,V_5)\in\cp^1\T\cp^2\T\cp^3\T\cp^4\T\cp^5$ 
such that $V_1\sub V_2\sub V_3\sub V_4\sub V_5$.
(Such a sequence is identified with its stabilizer $\{g\in G;gV_1=V_1,gV_2=V_2,gV_4=V_4,g\tV_4=\tV_4\}$.)

\subhead 3.21\endsubhead
Let
$$M'=\{(V_1,V_2,V_3,V_4,V_5)\in\cb;V_2=K_2,V_4=K_4\}.$$
If $(V_1,V_2,V_3,V_4,V_5)\in M$ then 

$AV_5\sub AV=K_4\sub V_5$, $AV_3\sub AK_4=A\ker(A^2)\sub\ker(A)=K_2\sub V_3$,

$AV_1\sub AK_2=0$, $AV_4=AK_4=K_2\sub V_4$, $AK_2=0$. 
\nl
We see that $(V_1,V_2,V_3,V_4,V_5)\in\cb_u$. Thus, 
$M'\sub\cb_u$. Note that $M'$ is isomorphic to $\PP^1\T\PP^1\T \PP^1$ hence it is a closed irreducible subvariety 
of $\cb_u$ of dimension $3$. Thus it is an irreducible component of $\cb_u$. It satisfies $J_{M'}=\{1,3,5\}$.
Hence $M'=M$.

\subhead 3.22\endsubhead
Let 
$$\align&T'=\{(V_1,V_2,V_3,V_4,V_5)\in\cb;V_1\sub K_2,K_4\sub V_5,A^2V_5=V_1,\\&
AV_4=V_2,V_2\sub AV_5\sub V_4\},\endalign$$
$$\align&Z=\{(V_1,V_2,V_4,V_5)\in\cp^1\T\cp^2\T\cp^4\T\cp^5;V_1\sub V_2,V_4\sub V_5,\\&
V_1\sub K_2,K_4\sub V_5,A^2V_5=V_1,AV_4=V_2,V_2\sub AV_5\sub V_4\}.\endalign$$
The map $\z:T'@>>>Z$, $(V_1,V_2,V_3,V_4,V_5)\m(V_1,V_2,V_4,V_5)$ is a $\PP^1$-bundle. Let 
$$Z'=\{(V_1,V_5)\in\cp^1\T\cp^5;V_1\sub K_2,K_4\sub V_5,A^2V_5=V_1\}.$$
Consider the map $\z':Z@>>>Z'$, $(V_1,V_2,V_4,V_5)\m(V_1,V_5)$.
Its fibre at $(V_1,V_5)$ can be identified with the projective line
$\{V_4; AV_5\sub V_4\sub V_5\}$ of the $2$-dimensional vector space $V_5/AV_5$ 
(note that $V_2$ is uniquely determined by $V_4$ via $V_2=AV_5$ and it automatically satisfies
$V_1\sub V_2\sub AV_5$ since $A^2V_5\sub AV_4\sub AV_5$); thus $\z'$ is a $\PP^1$-bundle. Note also that
$Z'$ can be identified  with the projective line $\{V_5; K_4\sub V_5\sub V\}$ of the $2$-dimensional vector 
space $V/K_4$. We see that $T'$ is a $\PP^1$-bundle over a $\PP^1$-bundle over a $\PP^1$-bundle. 
In particular, $T'$ is a closed smooth irreducible subvariety of dimension $3$ of $\cb$.
We show that $T'\sub\cb_u$. It is enough to show that if $(V_1,V_2,V_3,V_4,V_5)\in T'$, then $AV_i\sub V_i$ for 
$i\in I$.
For $i=5$ this follows from $AV_5\sub V_4\sub V_5$. For $i=4$ this follows from $AV_4=V_2\sub V_4$.
For $i=3$ this follows from $AV_3\sub AV_4=V_2\sub V_3$. For $i=2$ this follows from $AV_2\sub AV_4=V_2$.
For $i=1$ this follows from $AV_1\sub AK_2=0$.
Since $d_u=3$, we see that $T'\in\un\cb_u$. Hence $T'_i$ is defined for $i\in I$.

\subhead 3.23\endsubhead
Note that $T'_3=T'$. From the definitions we have
$$T'_1=\{(V_1,V_2,V_3,V_4,V_5)\in T'; V_2=K_2\},$$
$$T'_2=\{(V_1,V_2,V_3,V_4,V_5)\in T'; AV_3\sub V_1\},$$
$$T'_4=\{(V_1,V_2,V_3,V_4,V_5)\in T'; AV_5\sub V_3\},$$
$$T'_5=\{(V_1,V_2,V_3,V_4,V_5)\in T'; V_4=K_4\}.$$
We show:
$$T'_2=T'_4,\tag a$$
$$T'_1=T'_5.\tag b$$
Let $(V_1,V_2,V_3,V_4,V_5)\in T'$. 

If $AV_5\sub V_3$ then, since $\dim AV_5=3$, we have $AV_5=V_3$. Hence $AV_3=A^2V_5=V_1$. Thus $T'_4\sub T'_2$. 

If $AV_3\sub V_1$ then $A^2V_3\sub AV_1=0$ hence $V_3\sub K_4$. Now $A:K_4@>>>K_2$ is surjective and 
$V_1\sub K_2$ hence $\{x\in K_4;Ax\in V_1\}$ is a $3$-dimensional subspace of $K_4$ containing $V_3$ hence
$\{x\in K_4;Ax\in V_1\}=V_3$. Now $AV_5\sub\{x\in K_4;Ax\in V_1\}$ since $AV_5\sub K_4$ and $A^2V_5=V_1$.
Hence $AV_5\sub V_3$ (and even $AV_5=V_3$). Thus $T'_2\sub T'_4$.

Assume that $V_4=K_4$. Now $AV_4=V_2$, $AK_4=K_2$ hence $V_2=K_2$. Thus $T'_5\sub T'_1$.

Assume that $V_2=K_2$. We have $K_4=A\i(K_2)$, $V_4\sub A\i(V_2)=A\i(K_2)$ so that $V_4\sub K_4$. We see that
$T'_1\sub T'_5$.

This proves (a),(b). 
We set $T'_{24}=T'_2=T'_4$, $T'_{15}=T'_1=T'_5$. We have
$$T'_{24}=\{(V_1,V_2,V_3,V_4,V_5)\in\cb;V_1\sub K_2,K_4\sub V_5,AV_5=V_3, AV_3=V_1,AV_4=V_2\},$$
$$T'_{15}=\{(V_1,V_2,V_3,V_4,V_5)\in\cb;V_2=K_2,V_4=K_4,A^2V_5=V_1\}.$$
We set
$$\fZ:=T'_{24}\cap T'_{15}
=\{(V_1,V_2,V_3,V_4,V_5)\in\cb;AV_5=V_3, AV_3=V_1,V_2=K_2,V_4=K_4\}.$$
Now $\vt:T'_{24}@>>>\fZ$, $(V_1,V_2,V_3,V_4,V_5)\m(V_1,K_2,V_3,K_4,V_5)$ is a $\PP^1$-bundle; its fibre at 
$(V_1,K_2,V_3,K_4,V_5)$ is
$\{(V_2,V_4)\in\cp^2\T\cp^4;V_1\sub V_2\sub V_3, V_3\sub V_4\sub V_5, AV_4=V_2\}$ which is isomorphic to the 
projective line  $\{V_4;V_3\sub V_4\sub V_5\}$ of the $2$-dimensional vector space $V_5/V_3$.

Note that $\fZ$
is isomorphic to the projective line $\{V_5; K_4\sub V_5\sub V\}$ of the $2$-dimensional vector space $V/K_4$. 
Thus $T'_{24}$ is a $\PP^1$-bundle over $\PP^1$ and thus, having dimension $2$, is $\ne T'$. Now 
$$T'_{15}@>>>\{(V_1,V_2,V_4,V_5)\in\cp^1\T\cp^2\T\cp^4\T\cp^5;V_2=K_2,V_4=K_4,A^2V_5=V_1\}$$
is a $\PP^1$-bundle whose base is isomorphic to the projective line $\{V_5;K_4\sub V_5\sub V\}$ of the 
$2$-dimensional vector space $V/K_4$. Thus $T'_{15}$ is a $\PP^1$-bundle over 
$\PP^1$ and thus, having dimension $2$, 
is $\ne T'$. Also $T'_{15}$ is a union of $3$-lines (the fibres of the map above.)

We now see that $T'_i\ne T'$ for $i\in I-\{3\}$. It follows that $T'=T$, see 3.19.

Now any $3$-line in $T$ is the fibre of $\z:T'@>>>Z$ at some $(V_1,V_2,V_4,V_5)\in Z$. This fibre contains 
exactly one point in $T'_{24}$ namely the one defined by $V_3=AV_5$. This completes the proof 
of 3.19(a).

\subhead 3.24\endsubhead
In this subsection we assume that $G_{ad}$ is of type $A_2A_2A_1$. We can write $I=\{0,1,2,4,5\}$ where the 
notation is chosen so that $s_1s_2,s_4s_5$ have order $3$. We shall assume that $u\in G$ (see 
1.1) is such that the image of $u$ in $G_{ad}$ has a projection to each of the three factors a 
subregular element in that factor. We have $d_u=3$. Now 
$\cb_u$ has four irreducible components: $C_1,C_2,C_3,C_4$ where $J_{C_1}=\{1,0,4\}$, $J_{C_2}=\{2,0,4\}$, 
$J_{C_3}=\{2,0,5\}$, $J_{C_4}=\{1,0,5\}$. It is clear that there are well defined $Q^i\in\cp^i$, $i=1,2,4,5,$ 
such that $\cap_{i\in\{1,2,4,5\}}Q^i$ contains a Borel subgroup and
$$C_1=\{B\in\cb_u;B(2)=Q^2,B(5)=Q^5\},$$
$$C_2=\{B\in\cb_u;B(1)=Q^1,B(5)=Q^5\},$$
$$C_3=\{B\in\cb_u;B(1)=Q^1,B(4)=Q^4\}.$$
$$C_4=\{B\in\cb_u;B(2)=Q^2,B(4)=Q^4\},$$

\subhead 3.25\endsubhead
In this subsection we assume that $G_{ad}$ is of type $A_2A_1A_1$. We can write $I=\{0,1,3,5\}$ where 
the notation is chosen so that $s_0s_3$ has order $3$. We shall assume that $u\in G$ (see 1.1) is such 
that the image of $u$ in $G_{ad}$ has a projection to each of the three factors a subregular element in that factor. We have 
$d_u=3$. Now $\cb_u$ has two irreducible components: $C'_1,C'_2$ where $J_{C'_1}=\{1,3,5\}$, $J_{C'_2}=\{1,0,5\}$. 
It is clear that there are well defined parabolic subgroups
$P^0\in\cp^0$, $P^3\in\cp^3$ such that $P^0\cap P^3$ contains a Borel subgroup and
$$C'_1=\{B\in\cb_u;B(0)=P^0\},\qua C'_2=\{B\in\cb_u;B(3)=P^3\}.$$

\head 4. Euler characteristic computations\endhead
\subhead 4.1 \endsubhead
In this section we assume that $G_{ad}$ is of type $E_6$. We can write $I=\{0,1,2,3,4,5\}$ where 
the numbering is chosen so that $s_1s_2,s_2s_3,s_3s_4,s_4s_5,s_0s_3$ have order $3$. We shall 
assume that $u\in G$ (see 1.1) is such that $d_u=3$; this determines $u$ uniquely up to conjugacy. 
The Springer representation $\r_u$ of $W$ is a direct sum of two irreducible representations: one 
of dimension $30$ (appearing in the third symmetric power of the reflection representation) and 
one of dimension $15$ (appearing in the fifth symmetric power of the reflection representation).
Using 1.4(b) and the knowledge of $\r_u$ we see that there are exactly 

-two irreducible components $X$ of $\cb_u$ such that $J_X=\{3\}$ (we call them $\SS,\TT$);

-two irreducible components $X$ of $\cb_u$ such that $J_X=\{0,3\}$ (we call them $X(03),X'(03)$);

-one irreducible component $X$ of $\cb_u$ such that $J_X=\{0,1,4\}$ (we call it $X(014)$);

-one irreducible component $X$ of $\cb_u$ such that $J_X=\{0,2,4\}$ (we call it $X(024)$);

-one irreducible component $X$ of $\cb_u$ such that $J_X=\{0,2,5\}$ (we call it $X(025)$);

-one irreducible component $X$ of $\cb_u$ such that $J_X=\{0,1,5\}$ (we call it $X(015)$);

-one irreducible component $X$ of $\cb_u$ such that $J_X=\{1,3,5\}$ (we call it $X(135)$).

From 1.1(b) we see that if $K$ is one of 
$\{0,1,4\}$, $\{0,2,4\}$, $\{0,2,5\}$, $\{0,1,5\}$, $\{1,3,5\}$ then the exists a unique $P_K\in\cp_K$ such that 
$X(K)=\{B\in\cb;B\sub P_K\}$.

Now $u$ is induced from a unipotent element $u_1\in\bP$ where $P\in\cp_{01245}$ and $(u_1,\bP)$ is like
$(u,G)$ in 3.24; in particular we have $u\in P,u_1=\p_P(u)$. 
The four irreducible components of $\cb^{\bP}_{u_1}$ (see 3.24) give rise as in 1.5 to
four irreducible components of $\cb_u$ with the $J_X$ being preserved; these irreducible components of $\cb_u$
must be the same as $X(014)$, $X(024)$, $X(025)$, $X(015)$. Using 3.24, we see that there exist
$Q^i\in\cp^i$ (relative to $G$) where $i\in\{1,2,3,4,5\}$ such that
$$X(014)=\{B\in\cb;B(2)=Q^2,B(3)=Q^3,B(5)=Q^5\},$$
$$X(024)=\{B\in\cb;B(1)=Q^1,B(3)=Q^3,B(5)=Q^5\},$$
$$X(025)=\{B\in\cb;B(1)=Q^1,B(3)=Q^3,B(4)=Q^4\},$$
$$X(015)=\{B\in\cb;B(2)=Q^2,B(3)=Q^3,B(4)=Q^4\}.$$
(Note that these conditions determine $Q^i$ uniquely and that $Q^3=P$.)
Moreover, $\cap_{i\in\{1,2,3,4,5\}}Q^i$ contains a Borel subgroup.
Next we note that $u$ is induced from a unipotent element $u_2\in\bP'$ where $P'\in\cp_{0135}$ and 
$(u_2,\bP')$ is like $(u,G)$ in 3.25; in particular we have $u\in P',u_2=\p_{P'}(u)$. 
The two irreducible components of $\cb^{\bP'}_{u_2}$ (see 3.25) give rise as 
in 1.5 to two irreducible components of $\cb_u$ with the $J_X$ being preserved; these irreducible components of 
$\cb_u$ must be the same as $X(135)$, $X(015)$. Using 3.25 we see that there exist
$\tQ^i\in\cp^i$ (relative to $G$) where $i\in\{0,2,3,4\}$ such that
$$X(135)=\{B\in\cb;B(0)=\tQ^0,B(2)=\tQ^2,B(4)=\tQ^4\},$$
$$X(015)=\{B\in\cb;B(2)=\tQ^2,B(3)=\tQ^3,B(4)=\tQ^4\}.$$
(Note that these conditions determine $\tQ^i$ uniquely and that $\tQ^2$, $\tQ^4$ contain $P'$.)
Moreover, $\cap_{i\in\{0,2,3,4\}}\tQ^i$ contains a Borel subgroup.
Since we have 
$$\align&X(015)=\{B\in\cb;B(2)=Q^2,B(3)=Q^3,B(4)=Q^4\}\\&
=\{B\in\cb_u;B(2)=\tQ^2,B_3=\tQ^3,B(4)=\tQ^4\}\endalign$$
we see that we have $Q^2=\tQ^2$, $Q^3=\tQ^3$, $Q^4=\tQ^4$. We set $Q^0=\tQ^0$. Then we have
$$X(135)=\{B\in\cb;B(0)=Q^0,B(2)=Q^2,B(4)=Q^4\}.$$
From the previous argument we see also that $Q_0:=Q^1\cap Q^2\cap Q^3\cap Q^4\cap Q^5$ is a parabolic subgroup 
in $\cp_0$ and $Q_{15}:=Q^0\cap Q^2\cap Q^3\cap Q^4$ is a parabolic subgroup in $\cp_{15}$. Let 
$Q_{015}=Q^2\cap Q^3\cap Q^4$. We have $Q_{015}\in\cp_{015}$ and $Q_0\sub Q_{015}$, $Q_{15}\sub Q_{015}$. Since 
the adjoint group of $\bQ_{015}$ is of type $A_1A_1A_1$, the intersection $Q_0\cap Q_{15}$ is a Borel subgroup 
of $Q_{015}$. We see that $Q^0\cap Q^1\cap Q^2\cap Q^3\cap Q^4\cap Q^5$ is a Borel subgroup of $G$. We show:

(a) {\it The intersection $X(135)\cap X(024)$ consists of exactly one Borel subgroup $B_0$.}
\nl
Indeed the condition that $B\in\cb$ belongs to $X(135)\cap X(024)$ is the same as
$B(0)=Q^0,B(2)=Q^2,B(4)=Q^4,B(1)=Q^1,B(3)=Q^3,B(5)=Q^5$ that is,
$B\sub Q^0\cap Q^1\cap Q^2\cap Q^3\cap Q^4\cap Q^5$. It remains to use that the last intersection is a Borel
subgroup of $G$.

\subhead 4.2\endsubhead
Now $u$ is induced from a unipotent element $u'\in\bQ'$ where $Q'\in\cp_{12345}=\cp^0$ and 
$(u',Q')$ is like $(u,G)$ in 3.19; in particular, we have $u\in Q',u'=\p_{Q'}(u')$.
The irreducible components $M$ (resp. $T$)  of $\cb^{\bQ'}_{u'}$ as in 3.19 (with $G,u$ replaced by $\bQ',u')$)
give rise as in 1.5 to irreducible components of $\cb_u$ with the same $J_X$ which
therefore must be equal to $X(135)$ (resp. to $X'\in\{\SS,\TT\}$; we arrange notation so that $X'=\TT$).
From 3.19(a) we see that the following holds:

(a) {\it there are well defined parabolic subgroups $P'{}^2\in\cp^2$, $P'{}^4\in\cp^4$ such that
$X(135)=\{B\in\cb;B(0)=Q,B(2)=P'{}^2,B(4)=P'{}^4\}$;

we have $\TT\sub\{B\in\cb; B(0)=Q', B(1)\spa P'{}^2,B(5)\spa P'{}^4\}$;

we have $\TT_1=\TT_5=\{B\in\TT; B(0)=Q',B(2)=P'{}^2,B(4)=P'{}^4\}=X(135)\cap\TT$ (we denote it by $\TT_{15}$); 
this 
is a $\PP^1$-bundle over $\PP^1$ and is a union of $3$-lines;

we have $\TT_2=\TT_4$ (we denote it by $\TT_{24}$); it intersects any $3$-line in $\TT$ in exactly one point; 

the intersection $\TT_{15}\cap\TT_{24}$ is isomorphic to $\PP^1$;

there is a unique morphism $\vt:\TT_{24}@>>>\TT_{15}\cap\TT_{24}$ such that for any $B\in\TT_{24}$ we have 
$(B,\vt(B))\in\co_w$ for some $w\in W_{24}$; moreover, $\vt$ is a 
$\PP^1$-bundle.}
\nl
Since we have also $X(135)=\{B\in\cb;B(0)=Q^0,B(2)=Q^2,B(4)=Q^4\}$ we see that we must have
$Q'=Q^0,P'{}^2=Q^2,P'{}^4=Q^4$. Hence we have $\TT\sub\{B\in\cb; B(0)=Q^0, B(1)\spa Q^2,B(5)\spa Q^4\}$.

\subhead 4.3\endsubhead
Now $u$ is induced from a unipotent element $u''\in\bQ''$ where $Q''\in\cp_{0234}$ and $(u'',\bQ'')$ is like 
$(u,G)$ in 3.1 with $\a=0,\b=2,\g=4,\o=3$; in particular we have $u\in Q'',u''=\p_{Q''}(u)$.
The irreducible components $\hS$, $S_{\b\g'},S_{\b',\g}$ (resp. $S$)  of $\cb^{\bQ''}_{u''}$ 
defined as in 3.1 
(with $G,u$ replaced by $\bQ'',u'')$) give rise as in 1.5 to irreducible components of $\cb_u$ with the same 
$J_X$ which therefore must be equal to $X(024)$, $X(03)$, $X'(03)$ (resp. $X''\in\{\SS,\TT\}$); since
$\c(\SS)=\c(S)=12$ and $\c(\TT)=8$ (see 3.10(c), 3.19(b)), we must have $X''=\SS$.
We set $Q''=P^1\cap P^5$ where $P^1\in\cp^1,P^5\in\cp^5$. 
Let $P^0\ne\tP^0$ in $\cp^0$, $P^2\ne\tP^2$ in $\cp^2$, $P^4\ne\tP^4$ in $\cp^4$, $P^3\in\cp^3$ 
(with $\cp^i$ refering to $G$) be the inverse images under $\p_{Q''}$ of the parabolic subgroups with the same
names in $\bQ''$.

From 3.1(a) we see that the following holds:

(a) {\it We have}
$$X(024)=\{B\in\cb;B(1)=P^1,B(3)=P^3,B(5)=P^5\}$$
$$X(03)=\{B\in\cb;B(1)=P^1,B(2)=P^2,B(4)=\tP^4,B(5)=P^5\},$$
$$X'(03)=\{B\in\cb;B(1)=P^1,B(2)=\tP^2,B(4)=P^4,B(5)=P^5\},$$
$$\SS\sub\{B\in\cb;B(1)=P^1,B(5)=P^5\}.$$
Since
$$X(024)=\{B\in\cb;B(1)=Q^1,B(3)=Q^3,B(5)=Q^5\},$$
we see that $P^1=Q^1$, $P^3=Q^3$, $P^5=Q^5$. Thus,
$$X(03)=\{B\in\cb;B(1)=Q^1,B(2)=P^2,B(4)=\tP^4,B(5)=Q^5\},$$
$$X'(03)=\{B\in\cb;B(1)=Q^1,B(2)=\tP^2,B(4)=P^4,B(5)=Q^5\}.$$
From 3.1(e) we deduce

(b) {\it $X(024)\cap X(03)$ is a $0$-line;  $X(024)\cap X'(03)$ is a $0$-line.}

\subhead 4.4\endsubhead
Now $u$ is induced from a unipotent element ${}'u\in{}'\bQ$ where ${}'Q\in\cp_{01234}=\cp^5$ and 
$({}'u,{}'\bQ)$ is like $(u,G)$ in 3.15 with $\d=1,\a=2,\b=0,\g=4,\o=3$; in particular we have 
$u\in{}'Q,{}'u=\p_{{}'Q}(u)$. The irreducible components $\tS$, $E$ of $\cb^{{}'\bQ}_{{}'u}$ 
defined in 3.15 
(with $G,u$ replaced by ${}'\bQ,{}'u)$) give rise as in 1.5 to the irreducible components $\SS$ and
$X(014)$ of $\cb_u$. From 3.16(a),(b),(c) 
we deduce that there is a well defined $\hP^2\in\cp^2$ such that
$$X(014)=\{B\in\cb;B(2)=\hP^2,B(3)=Q^3,B(5)={}'Q\},$$
$$\hP^2\n\{P^2,\tP^2\},$$
$$\SS_1=\{B\in\SS;B(2)=\hP^2.\}$$
Since we have also
$$X(014)=\{B\in\cb;B(2)=Q^2,B(3)=Q^3,B(5)=Q^5\},$$
we see that $\hP^2=Q^2$, ${}'Q=Q^5$. Thus we have
$$Q^2\n\{P^2,\tP^2\},$$
$$\SS_1=\{B\in\SS;B(2)=Q^2.\}$$

\subhead 4.5\endsubhead
Now $u$ is induced from a unipotent element ${}''u\in{}''\bQ$ where ${}''Q\in\cp_{02345}=\cp^1$ and 
$({}''u,{}''\bQ)$ is like $(u,G)$ in 3.15
 with $\d=5,\a=4,\o=3,\b=2,\g=0$; in particular we have 
$u\in{}''Q,{}''u=\p_{{}''Q}(u)$. The irreducible components $\tS$, $E$  of $\cb^{{}''\bQ}_{{}''u}$ defined in 3.15
(with $G,u$ replaced by ${}''\bQ,{}''u)$) give rise as in 1.5
 to the irreducible components $\SS$ and $X(025)$ of
$\cb_u$. From 3.16(a),(b),(c) we deduce as in 4.4 that:
$$Q^4\n\{P^4,\tP^4\},$$
$$\SS_5=\{B\in\SS;B(4)=Q^4.\}$$

\subhead 4.6\endsubhead
We show: 

(a) $X(135)\cap X(03)=\emp,X(135)\cap X'(03)=\emp$.
\nl
Recall that

$X(135)=\{B\in\cb;B(0)=Q^0,B(2)=Q^2,B(4)=Q^4\}$, 

$X(03)=\{B\in\cb;B(1)=Q^1,B(2)=P^2,B(4)=\tP^4,B(5)=Q^5\},$

$X'(03)=\{B\in\cb;B(1)=Q^1,B(2)=\tP^2,B(4)=P^4,B(5)=Q^5\}.$
\nl
If $B\in X(135)\cap X(03)$ then $B(2)=Q^2=P^2$. This contradicts $Q^2\ne P^2$.
Similarly if $B\in X(135)\cap X'(03)$ then $B(2)=Q^2=\tP^2$. This contradicts $Q^2\ne\tP^2$. This proves (a).

Using $\TT_{15}=\TT\cap X(135)$ we see that (a) implies:

(b) $\TT_{15}\cap X(03)=\emp,\TT_{15}\cap X'(03)=\emp$.

\subhead 4.7\endsubhead
We show:

(a) {\it The intersection of a $0$-line in $\cb$ with $\TT$ is either a point or is empty.}
\nl
Assume that this intersection contains two distinct points $B',B''$. We have 

$(B',B'')\in \co_{s_0}$.
\nl
Since $\TT\sub Q^0$, we have $B'\in Q^0,B''\in Q^0$. Since 

$Q^0\in\cp_{12345}$,
\nl
 we have $(B',B'')\in \co_w$ where $w\in W_{12345}$. This contradicts $w=s_0$; (a) is proved.

We show:

(b) {\it The intersection $\TT\cap X(024)\cap X(03)$ is either a point or is empty.
The intersection $\TT\cap X(024)\cap X'(03)$ is either a point or is empty.}
\nl
This follows from (a) since $X(024)\cap X(03)$ is a $0$-line and $X(024)\cap X'(03)$ is a $0$-line, see 4.3(b).

\subhead 4.8\endsubhead
We write $X(135)\cap X(024)=\{B_0\}$ as in 4.1(a). We show:

(a) $\TT_{15}\cap\TT_{24}\cap\TT_0\sub\{B_0\}$.
\nl
Let $B\in\TT_{24}\cap\TT_0$. Let $P\in\cp^{024}$ be such that $B\sub P$. Since
$B\in\TT_0\cap\TT_2\cap\TT_4$, we see from  1.3(a) that $X=\cb^P$ is an irreducible component 
of $\cb_u$ and that
$J_X=\{0,2,4\}$. Thus we have $X=X(024)$. We see that $B\in X(024)$.
Thus $\TT_{24}\cap\TT_0\sub X(024)$. Recall that $\TT_{15}=X(135)\cap\TT$. Thus $\TT_{15}\sub X(135)$ and
$\TT_{15}\cap\TT_{24}\cap\TT_0\sub X(024)\cap X(135)=\{B_0\}$; now (a) follows.

\subhead 4.9\endsubhead
We show:

(a) {\it The morphism $\vt:\TT_{24}@>>>\TT_{24}\cap\TT_{15}$, see 3.19, restricts to a morphism
$\vt_0:\TT_{24}\cap\TT_0@>>>\TT_{24}\cap\TT_{15}\cap\TT_0$. Moreover, $\vt_0$ is a 
$\PP^1$-bundle.}
\nl
Let $B\in\TT_{24}\cap\TT_0$ and let $B'=\vt(B)\in\TT_{24}\cap\TT_{15}$. 
As in the proof of 4.8(a) we have $B\in X(024)$. Recall that
$(B,B')\in\co_w$ for some $w\in W_{24}$. 
Since $B\in X(024)$ it follows that $B'\in X(024)$. Thus $B'\in\TT_{24}\cap\TT_{15}\cap\TT_0$. Thus, $\vt_0$ is
well defined. 
Conversely, let $B'\in \TT_{24}\cap\TT_{15}\cap\TT_0$ and let $B\in\vt\i(B')$. Again we have $B'\in X(024)$
and $(B,B')\in\co_w$ for some $w\in W_{24}$. It follows that 
$B\in X(024)$ so that $B\in\TT_{24}\cap\TT_0$. We see that $\vt_0$ is a $\PP^1$-bundle. This proves (a).

We show:

(b) {\it If $B_0\in\TT$ (hence $B_0\in\TT_{24}\cap\TT_{15}\cap\TT_0$) then 
$\TT_{24}\cap\TT_0$ is a projective line whose intersection with any $3$-line in $\TT$ has exactly one point or
is empty. If $B_0\n\TT$, then $\TT_{24}\cap\TT_0=\emp$.}
\nl
Using (a), we see that it is enough to show the following statement. 

{\it If $B,B'$ are  Borel subgroups in $\TT_{24}\cap\TT_0$ such that $B,B'$ are on the same $3$-line in $\cb$ then
$B=B'$.}
\nl
Since $\vt(B)=\vt(B')=B_0$, we have $(B,B_0)\in\co_w$, $(B',B_0)\in\co_{w'}$, for some $w,w'$ in $W_{24}$.
 It follows that $(B,B')\in\co_{y}$ for some $y\in W_{24}$. Since $B,B'$ are on the same $3$-line we have also 
$(B,B')\in\co_z$ for some $z\in W_3$.  This forces $y=z=1$. This proves (b).

\subhead 4.10\endsubhead
We identify $Y\sub S$ with a subvariety of $\SS$. We show: 

(a) {\it If $B_0\n\TT$ then $\SS\cap\TT=\emp$. If $B_0\in\TT$ then $\SS\cap\TT$ is the union of all $3$-lines 
which intersect $\{B\in\TT_{24}\cap T_0;\vt(B)=B_0\}\cap Y$.}
\nl
For $K\sub I$ we denote by $w(K)$ an element of $W_K$. Let $B\in\SS\cap\TT$. The $3$-line 
through $B$ must intersect $Y$ in a unique point $B_1$ and it must intersect $\TT_{24}$ in a 
unique point $B_2$. Let $B_3=\vt(B_2)\in\TT_{24}\cap\TT_{15}\sub X(135)$. We have 
$(B_3,B_0)\in\co_{w(135)}$. We have $Y\sub X(024)$, hence $(B_1,B_0)\in\co_{w(024)}$. We have 
$(B,B_1)\in\co_{w(3)}$, hence $(B,B_0)\in\co_{w(3)w(024)}$.

We have $(B,B_2)\in\co_{w'(3)}$, $(B_2,B_3)\in\co_{w(24)}$, $(B_3,B_0)\in\co_{w(135)}$.
If $w(24)\ne1$ it follows that $(B,B_0)\in\co_{w'(3)w(24)w(135)}$ hence $w'(3)w(24)w(135)=
w(3)w(024)$ so that
$w(135)=1,w'(3)=w(3),w(024)=w(24)$ and $B_3=B_0$, $B_1=B_2\in\TT_{24}\cap Y$, $(B_1,B_0)\in\co_{w(24)}$. 
Since $B_3\in\TT$ and $B_3=B_0$ we have also $B_0\in\TT$ (hence $B_0\in\TT_{24}\cap\TT_{15}\cap\TT_0$) and
$B_1\in\{B'\in\TT_{24}\cap T_0;\vt(B)=B_0\}\cap Y$.

If $w_{24}=1$ we must have $B_2=B_3$, $(B_3,B_0)\in\co_{w_1(3)}$. Since $B_3\in\TT$ and $\TT$ is a union of
$3$-lines it follows that $B_0\in\TT$. From $(B,B_2)\in\co_{w'(3)},(B_2,B_0)\in\co_{w_1(3)}$ we obtain
$(B,B_0)\in\co_{w_2(3)}$.

Now (a) follows (we have used 4.9(b)).

\subhead 4.11\endsubhead
We show: 

(a) {\it Assume that $B_0\in\TT$. Let $U=\{B\in\TT_{24}\cap T_0;\vt(B)=B_0\}$ (a projective line).
Then either $U\cap Y\cong\PP^1$ or $U\cap Y$ consists of two points or $U\cap Y$ consists of one point.} 
\nl
Recall that $Y$ can be viewed as the variety
$$\{(z_0,z'_0,z_2,z'_2),(y,y'))\in\PP^3\T \PP^1;z_0z'_0+z_2z'_2=0, z_2y'+z'_2y=0\}$$
or setting $z_0=ac,z'_0=bd,z_2=-ad,z'_2=bc$, as the variety
$$\{((a,b),(c,d),(y,y'))\in\PP^1\T\PP^1\T\PP^1;-ady'+bcy=0\}.$$
Now $U$ can be viewed as a subvariety of $\PP^1\T\PP^1\T\PP^1$ of the form
$$\{((a,b),(c,d),(y,y'))\in\PP^1\T\PP^1\T\PP^1; c=m_1a+m_2b, d=m'_1a+m'_2b, y=y_0,y'=y_0\}$$
where $(y_0,y'_0)\in\PP^1$ is fixed and $(m_1,m_2,m'_1,m'_2)\in\kk^4$ is fixed such that $m_1m'_2-m_2m'_1\ne0$.
Then $U\cap Y$ becomes
$$\align&
\{((a,b),(c,d),(y,y'))\in\PP^1\T\PP^1\T\PP^1;\\&
-a(m'_1a+m'_2b)y'_0+b(m_1a+m_2b)y_0=0, y=y_0,y'=y_0\}
\\&\cong
\{((a,b)\in\PP^1;-a(m'_1a+m'_2b)y'_0+b(m_1a+m_2b)y_0=0\}\\&=
\{((a,b)\in\PP^1;-m'_1y'_0a^2+(-m'_2y'_0+m_1y_0)ab+m_2y_0b^2=0\}.\endalign$$
It follows that, if each of $-m'_1y'_0,-m'_2y'_0+m_1y_0,m_2y_0$ is zero then $U\cap Y\cong\PP^1$; otherwise,
$U\cap Y$ consists of one or two points. This proves (a).

\subhead 4.12\endsubhead
We set $\SS^\di=\SS-(\SS_0\cup\SS_1\cup\SS_2\cup\SS_4\cup\SS_5)$. This is an open dense subset 
of $\SS$. 
We have the following result.

\proclaim{Proposition 4.13} We have $\c(\SS^\di)=1$.
\endproclaim
Under the isomorphism $S@>\si>>\SS$, $\SS_5$ (resp. $\SS_1$) corresponds to a closed subsets $R_5$ (resp. $R_1$) 
of $S$ and $\SS^\di$ corresponds to the open dense subset $S^\di:=S-(S_0\cup R_1\cup S_2\cup S_4\cup R_5)$ of $S$. It is enough to prove that 
$$\c(S^\di)=1.\tag a$$ 
The proof is given in 4.14-4.16.  

\subhead 4.14\endsubhead 
To prove that $\c(S^\di)=1$ we can identify our $S$ with $X\sub\cb$ as in 3.12(c) in such a way
that $S_0=X_\a,S_2=X_\b,S_4=X_\g$ and, if 

$A:V@>>>V$, $L_1,L'_1,{}^!L_4,{}^!L'_4,L_4^!,L'_4{}^!,I_2,K_6$
\nl
 are as in 3.2, 3.4, then $R_5$ is identified with 
$$\{(V_1,V_2,V_4,\tV_4)\in X;V_4=\hV_4^!\},$$
 $R_1$ is identified with 
$$\{(V_1,V_2,V_4,\tV_4)\in X;\tV_4={}^!\hV_4\},$$
 where $\hV_4^!$ is an isotropic $4$-space in 
$V$ in the same family as $L_4^!,L'_4{}^!$ but distinct from each of them and ${}^!\hV_4^!$ is 
an isotropic $4$-space in $V$ in the same family as ${}^!L_4,{}^!L'_4$ but distinct from each 
of them; moreover we have $I_2\sub{}^!\hV_4$,
$I_2\sub\hV_4^!$. Recall that $X_\a=X_\a^*\cup Y$, $X_\b=X_\b^*\cup Y$, $X_\g=X_\g^*\cup Y$, where 
$$\align&X_\a^*=\{(V_1,V_2,V_4,\tV_4)\in\cb; V_1\sub I_2,V_4={}^!L_4,\tV_4=L'_4{}^!\}\\&\sqc
\{(V_1,V_2,V_4,\tV_4)\in\cb; V_1\sub I_2,V_4={}^!L'_4,\tV_4=L_4^!\},\endalign$$
$$\align&X_\b^*=\{(V_1,V_2,V_4,\tV_4)\in\cb;I_2\sub V_4,V_1=L_1,\tV_4=L_4^!\}\\&\sqc
\{(V_1,V_2,V_4,\tV_4)\in\cb;I_2\sub V_4,V_1=L'_1,\tV_4=L'_4{}^!\},\endalign$$
$$\align&X_\g^*=\{(V_1,V_2,V_4,\tV_4)\in\cb;I_2\sub\tV_4,V_1=L_1,V_4={}^!L_4\}\\&
\sqc\{(V_1,V_2,V_4,\tV_4)\in\cb;I_2\sub\tV_4,V_1=L'_1,V_4={}^!L'_4\},\endalign$$
$$Y=\{(V_1,V_2,V_4,\tV_4)\in X;V_2=I_2\}.$$

\subhead 4.15\endsubhead 
As in 3.10 let $Z$ be the variety of all pairs $V_1,V_3$ of isotropic subspaces of $V$ of 
dimension $1$ and $3$ respectively such that $V_1\sub I_2\sub V_3$ and $AV_3\sub V_1$ and let 
$\fU$ be the space of isotropic $3$-spaces in $V$ that contain $I_2$ hence are contained in
$K_6$; we have $\fU\cong\PP^1\T\PP^1$. From 3.10(d) we have $\c(Z)=6$.

Let $Z_1=\{(V_1,V_3)\in Z; V_3\sub{}^!\hV_4\}$, $Z_5=\{(V_1,V_3)\in Z; V_3\sub\hV_4^!\}$.
Let $\fU_1=\{V_3\in\fU;V_3\sub{}^!\hV_4\}$, $\fU_5=\{V_3\in\fU;V_3\sub\hV_4^!\}$.
Note that $\fU_1\cong\PP^1$, $\fU^5\cong\PP^1$.
By the argument in the proof of 3.10(d), the map $Z@>>>\fU$, $(V_1,V_3)\m V_3$, restricts to 
an isomorphism $Z_1@>\si>>\fU_1$ and to an isomorphism $Z_2@>\si>>\fU_2$ (note that 
$\fU_1,\fU_2$ do not contain the exceptional points
${}^!L_4\cap L'_4{}^!$ and ${}^!L'_4\cap L_4^!$). We see that $Z_1\cong\PP^1$, $Z_5\cong\PP^1$.
It follows that
$$\c(Z_1)=\c(Z_5)=2.$$
We have $Z_1\cap Z_5=\{(V_1,V_3)\in Z; V_3\sub{}^!\hV_4\cap\hV_4^!\}$.
Since $I_2\sub {}^!\hV_4\cap\hV_4^!$, the intersection ${}^!\hV_4\cap\hV_4^!$ is a $3$-dimensional
isotropic subspace; we see that $Z_1\cap Z_5$ consists of a single element $(V_1,V_3)$ where
$V_3={}^!\hV_4\cap\hV_4^!$ and $V_1$ is uniquely determined by $V_3$. In particular,
$$\c(Z_1\cap Z_5)=1.$$
We have 
$$\c(Z_1\cup Z_5)=\c(Z_1)+\c(Z_5)-\c(Z_1\cap Z_5)=2+2-1=3.$$
Let
$$Z_\a^*=\{(V_1,V_3)\in Z;V_1\sub I_2,V_3={}^!L_4\cap L'_4{}^!\}\sqc
\{(V_1,V_3)\in Z;V_1\sub I_2,V_3={}^!L'_4\cap L_4^!\},$$
$$Z_\b^*
=\{(V_1,V_3)\in Z;V_1=L_1,V_3\sub L_4^!\}\sqc\{(V_1,V_3)\in Z;V_1=L'_1,V_3\sub L'_4{}^!\},$$
$$Z_\g^*
=\{(V_1,V_3)\in Z;V_1=L_1,V_3\sub {}^!L_4\}\sqc\{(V_1,V_3)\in Z;V_1=L'_1,V_3\sub {}^!L'_4\},$$
We have
$$Z_\a^*\cap Z_\b^*=\{(L_1,{}^!L'_4\cap L_4^!)\},\sqc\{(L'_1, {}^!L_4\cap L'_4{}^!)\},$$
$$Z_\a^*\cap Z_\g^*=\{(L_1,{}^!L_4\cap L'_4{}^!)\}\sqc\{(L'_1,{}^!L'_4\cap L_4^!)\},$$
$$Z_\b^*\cap Z_\g^*=\{(L_1,{}^!L_4\cap L_4^!\}\sqc\{(L'_1,{}^!L'_4\cap L'_4{}^!\},$$
$$Z_\a^*\cap Z_\b^*\cap Z_\g^*=\emp.$$
Thus, each of $Z_\a^*,Z_\b^*,Z_\g^*$ is a disjoint union of two copies of $\PP^1$; each of
$Z_\a^*\cap Z_\b^*$, $Z_\a^*\cap Z_\g^*$, $Z_\b^*\cap Z_\g^*$ consists of two points. We see 
that
$$\c(Z_\a^*)=\c(,Z_\b^*)=\c(Z_\g^*)=4,$$
$$\c(Z_\a^*\cap Z_\b^*)=\c(Z_\a^*\cap Z_\g^*)=\c(Z_\b^*\cap Z_\g^*)=2$$.
$$\c(Z_\a^*\cap Z_\b^*\cap Z_\g^*)=0.$$
We have
$$\align&\c(Z_\a^*\cup Z_\b^*\cup Z_\g^*)=\c(Z_\a^*)+\c(,Z_\b^*)+\c(Z_\g^*)\\&
-\c(Z_\a^*\cap Z_\b^*)-\c(Z_\a^*\cap Z_\g^*)
-\c(Z_\b^*\cap Z_\g^*)+\c(Z_\a^*\cap Z_\b^*\cap Z_\g^*)\\&=4+4+4-2-2-2+0=6.\endalign$$
We have
$$Z_\a^*\cap Z_1=\emp, Z_\g^*\cap Z_1=\emp,$$    
$$Z_\b^*\cap Z_1=\{(L'_1,{}^!\hV_4\cap L'_4{}^!)\}\sqc\{(L_1,{}^!\hV_4\cap L_4^!)\},$$
$$Z_\a^*\cap Z_5=\emp, Z_\b^*\cap Z_5=\emp,$$    
$$Z_\g^*=\{(L_1,{}^!L_4\cap\hV_4^!\}\sqc\{(L'_1,{}^!L'_4\cap\hV_4^!\}.$$
Thus 
$$(Z_\a^*\cup Z_\b^*\cup Z_\g^*)\cap Z_1
=\{(L'_1,{}^!\hV_4\cap L'_4{}^!)\}\sqc\{(L_1,{}^!\hV_4\cap L_4^!)\},$$
$$(Z_\a^*\cup Z_\b^*\cup Z_\g^*)\cap Z_5
=\{(L_1,{}^!L_4\cap\hV_4^!\}\sqc\{(L'_1,{}^!L'_4\cap\hV_4^!\}.$$
Hence
$$(Z_\a^*\cup Z_\b^*\cup Z_\g^*)\cap (Z_1\cap Z_5)=\emp.$$
We have 
$$\c((Z_\a^*\cup Z_\b^*\cup Z_\g^*)\cap Z_1)=\c((Z_\a^*\cup Z_\b^*\cup Z_\g^*)\cap Z_5)=2,$$
$$\c(Z_\a^*\cup Z_\b^*\cup Z_\g^*)\cap (Z_1\cap Z_5))=0$$
hence 
$$\align&\c((Z_\a^*\cup Z_\b^*\cup Z_\g^*)\cap (Z_1\cup Z_5))\\&=
\c((Z_\a^*\cup Z_\b^*\cup Z_\g^*)\cap Z_1)+\c((Z_\a^*\cup Z_\b^*\cup Z_\g^*)\cap Z_5)\\&-
\c(Z_\a^*\cup Z_\b^*\cup Z_\g^*)\cap (Z_1\cap Z_5))=2+2-0=4.\endalign$$
We have
$$\align&\c((Z_\a^*\cup Z_\b^*\cup Z_\g^*)\cup(Z_1\cup Z_5))\\&=
\c(Z_\a^*\cup Z_\b^*\cup Z_\g^*)+\c(Z_1\cup Z_5)\\&
-\c((Z_\a^*\cup Z_\b^*\cup Z_\g^*)\cap (Z_1\cup Z_5))=6+3-4=5.\endalign$$

\subhead 4.16\endsubhead
We have a partition $Z=Z'\cup Z''$ where $Z'=Z_\a^*\cup Z_\b^*\cup Z_\g^*\cup Z_1\cup Z_5$ and 
$Z''=Z-Z'$.
Recall that $\z:(V_1,V_2,V_4,\tV_4)\m(V_1,V_4\cap\tV_4)$ makes $X$ into a $\PP^1$-bundle over 
$Z$ and that  
$$Y=\{(V_1,V_2,V_4,\tV_4)\in X;V_2=I_2\}.$$
Recall that $\z$ restricts to an isomorphism $Y@>\si>>Z$. Hence setting $Y'=\z\i(Z')\cap Y$,
$Y''=\z\i(Z'')\cap Y$, we have a partition $Y=Y'\cup Y''$ and $\z$ restricts to isomorphisms
$Y'@>\si>>Z'$, $Y''@>\si>>Z''$. We have 
$$\c(Y'')=\c(Z'')=\c(Z)-\c(Z')=6-5=1.$$ 
We have
$$\align&X_\a^*=\z\i(Z_\a^*),\qua X_\b^*=\z\i(Z_\b^*),\qua X_\g^*=\z\i(Z_\g^*),\\&
\qua R_1=\z\i(Z_1),\qua R_5=\z\i(Z_5).\endalign$$
Hence 
$$X_\a^*\cup X_\b^*\cup X_\g^*\cup R_1\cup R_5=\z\i(Z'),$$
so that
$$\c(X_\a^*\cup X_\b^*\cup X_\g^*\cup R_1\cup R_5)=\c(\z\i(Z'))=2\c(Z')=10.$$
Recall from 3.1(c) that
$$X_\a\cup X_\b\cup X_\g\cup R_1\cup R_5=X^*_\a\cup X^*_\b\cup X^*_\g\cup R_1\cup R_5\cup Y$$
hence we have a partition
$$X_\a\cup X_\b\cup X_\g\cup R_1\cup R_5
=(X^*_\a\cup X^*_\b\cup X^*_\g\cup R_1\cup R_5)\cup Y'',$$
so that
$$\c(X_\a\cup X_\b\cup X_\g\cup R_1\cup R_5)=\c(X^*_\a\cup X^*_\b\cup X^*_\g\cup R_1\cup R_5)
+\c(Y'')=10+1=11.$$
We have $\c(X)=2\c(Z)=12$ (see 3.10(c),(d)) hence
$$\c(X-(X_\a\cup X_\b\cup X_\g\cup R_1\cup R_5))=\c(X)-\c(X_\a\cup X_\b\cup X_\g\cup R_1\cup R_5)
=12-11=1.$$
This proves 4.13(a). Proposition 4.13 is proved.

\subhead 4.17\endsubhead
We show:

(a) {\it Let $\cl$ be a $3$-line in $\TT$. Then either $\cl\sub\TT_0$ or $\cl\cap\TT_0$ is 
exactly one point.}
\nl
Let $B\in\cl$. Let $P\in\cp_{03}$ be such that $B\sub P$. We have $u\in P$. Let $u'=\p_P(u)$. 
Let $\bar\cl=\{\p_P(B');B'\in\cl\}$. Note that $\bar\cl$ is a line contained in 
$\cb^{\bP}_{u'}$. It follows that $u'$ is not regular unipotent in $\bP$. Since $\bP_{ad}$ is 
of type $A_2$, $u'$ must be subregular or $1$. If $u'=1$ then $\cb^{\bP}_{u'}=\cb^{\bP}$ so 
that for any $B'\in\cl$ the $0$-line through $B'$ is contained in $\cb^P$ hence is contained in
$\cb_u$; thus $\cl\sub\TT_0$. If $u'$ is subregular then 
$\cb^{\bP}_{u'}$ is the union of two projective lines which intersect in a single point. We see
that there is a unique $B'\in\cl$ such that the $0$-line through $B'$ is contained in $\cb_u$; 
thus $B'\in\TT_0$. This proves (a).

\subhead 4.18\endsubhead
We set $\TT^\di=\TT-(\TT_0\cup\TT_1\cup\TT_2\cup\TT_4\cup\TT_5)$. This is an open dense subset 
of $\TT$. We have the following result.

\proclaim{Proposition 4.19}If $B_0\n\TT$ then $\c(\TT^\di)=0$. If $B_0\in\TT$ then 
$\c(\TT^\di)=1$.
\endproclaim
Let $\fV$ be the variety of all $3$-lines in $\TT$. We have a partition
$\fV=\fV_{15}\sqc\hat\fV\sqc\fV'\sqc\fV''$ where

$\fV_{15}$ is the variety of all $3$-lines $\cl$ such that $\cl\sub\TT_{15}$;

$\hat\fV$ is the variety of all $3$-lines $\cl$ such that $\cl\sub\TT-\TT_{15}$, $\cl\sub\TT_0$;

$\fV'$ is the variety of all $3$-lines $\cl$ such that $\cl\sub\TT-\TT_{15}$, $\cl\cap\TT_0$ is
exactly one point and that point is in $\TT_{24}$;

$\fV''$ is the variety of all $3$-lines $\cl$ such that $\cl\sub\TT-\TT_{15}$, $\cl\cap\TT_0$ 
is exactly one point and that point is not in $\TT_{24}$.
\nl
The fact that this partition is well defined follows from 3.19(a) and the fact that $\TT_{15}$ 
is a union of $3$-lines. Let $\x:\TT@>>>\fV$ be the map which associates to $B\in\TT$ the 
$3$-line containing $B$.

The inverse images of $\fV_{15},\hat\fV,\fV',\fV'$ under $\x$ are denoted by 
$\TT_{15},\hat\TT,\TT',\TT''$. (This agrees with our earlier definition of $\TT_{15}$.) We have
$\TT=\TT_{15}\sqc\hat\TT\sqc\TT'\sqc\TT''$. We set $\TT'{}^\di=\TT'\cap\TT^\di$, 
$\TT''{}^\di=\TT''\cap\TT^\di$. Clearly, we have $\TT^\di=\TT'{}^\di\sqc\TT''{}^\di$. 
Let $\TT''{}^\di@>>>\fV''$ be the restriction of $\x$. This is a fibration with each fibre 
being a projective line from which two points have been removed: one in $\TT_{24}$ and one in 
$\TT_0$. It follows that $\c(\TT''{}^\di)=0$. Let $\TT'{}^\di@>>>\fV'$ be the restriction of 
$\x$. This is a fibration with each fibre being a projective line from which one point has been
removed: the one in $\TT_{24}$ (which is also in $\TT_0$). It follows that 
$\c(\TT'{}^\di)=\c(\fV')$. We have 
$$\c(\TT^\di)=\c(\TT''{}^\di)+\c(\TT'{}^\di)=0+\c(\fV')=\c(\fV').$$
If $B_0\n\TT$ then $\TT_{24}\cap\TT_0=\emp$ (see 4.9(b)) hence $\fV'=\emp$ and 
$\c(\TT^\di)=\c(\fV')=0$. Now assume that $B_0\in\TT$. Then by 4.9(b), $\fV'$ is isomorphic to
$\TT_{24}\cap\TT_0-\TT_{24}\cap\TT_{15}\cap\TT_0$ and this is a projective line minus a point. 
Thus, in this case, $\c(\TT^\di)=\c(\fV')=1$.

\subhead 4.20\endsubhead
We show:

(a) {\it If $B_0\n\TT$ then $\SS\cap\TT=\emp$; hence $\SS^\di\cap\TT^\di=\emp$ and 
$\c(\SS^\di\cap\TT^\di)=0$.}

(b) {\it If $B_0\in\TT\cap\SS$ then $\c(\SS^\di\cap\TT^\di)$ is $1$ or $0$.}

(c) {\it If $B_0\in\TT,B_0\n\SS$ then $\c(\SS^\di\cap\TT^\di)$ is $2$ or $1$.}
\nl
Now (a) follows from 4.10(a). Next we assume that $B\in\TT$. As in 4.11(a), let
$U=\{B\in\TT_{24}\cap T_0;\vt(B)=B_0\}$. We have a fibration 
$\SS^\di\cap\TT^\di@>>>(U\cap Y)-\{B_0\}$. This associates to $B$ the intersection of the 
$3$-line through $B$ with $\TT_{24}$, see 4.2(a); each fibre is a projective line with a point 
removed. It follows that $\c(\SS^\di\cap\TT^\di)=\c((U\cap Y)-\{B_0\})$. This equals 
$\c(U\cap Y)$ if $B_0\n Y$ (that is if $B_0\n\SS$) and it equals $\c(U\cap Y)-1$ if $B_0\in Y$ 
(that is if $B_0\in\SS$). It remains to use 4.11(a).

\subhead 4.21\endsubhead
We show:

(a) {\it If $B_0\n\TT$ then  $\c(\SS^\di\cup\TT^\di)=1$.}

(b) {\it If $B_0\in\TT$ then $\c(\SS^\di\cup\TT^\di)\in\{0,1,2\}$.}
\nl
We have $\c(\SS^\di\cup\TT^\di)=\c(\SS^\di)+\c(\TT^\di)-\c(\SS^\di\cap\TT^\di)$. It remains to
use 4.13, 4.19 and 4.20(a),(b),(c).

\subhead 4.22\endsubhead
Let $\cb_u^\di=\cb_u-(\cb_{u,0}\cup\cb_{u,1}\cup\cb_{u,2}\cup\cb_{u,4}\cup\cb_{u,5})$. We show:

(a) $\cb_u^\di=\SS^\di\cup\TT^\di$.
\nl
The inclusion $\SS^\di\cup\TT^\di\sub\cb_u^\di$ is obvious. Conversely, let $B\in\cb_u^\di$. It
is enough to show that $B\in\SS\cup\TT$. We have $B\in X$ for some $X\in\un\cb$. If $X$ is 
$\SS$ or $\TT$, we are done. Thus we can assume that $X\n\{\SS,\TT\}$. We have $J_X\ne\{3\}$ 
and since $J_X\ne\emp$ we have $i\in J_X$ for some $i\in I-\{3\}$ and, in particular, 
$X\sub\cb_{u,0}\cup\cb_{u,1}\cup\cb_{u,2}\cup\cb_{u,4}\cup\cb_{u,5}$. Since $B\in X$, we have 
$B\in\cb_{u,0}\cup\cb_{u,1}\cup\cb_{u,2}\cup\cb_{u,4}\cup\cb_{u,5}$. This contradicts
$B\in\cb_u^\di$. This proves (a). 

From (a) and 4.21(a),(b) we deduce

(b) {\it We have $\c(\cb_u^\di)\in\{0,1,2\}$.}

\subhead 4.23\endsubhead
Let $\cp_u^{reg}$ be the set of all $P\in\cp_{1245}$ such that $u\in P$ and $\p_P(u)$ is a 
regular unipotent element of $\bP$. We define a map $\ph:\cp_u^{reg}@>>>\cb_u$ by $P\m B$ where
$B$ is the unique Borel subgroup such that $u\in B\sub P$. Let 
$\cb_u^{reg}=\cb_u-(\cb_{u,1}\cup\cb_{u,2}\cup\cb_{u,4}\cup\cb_{u,5})$. We show:

(a) {\it $\ph$ defines an isomorphism $\ph_0:\cp_u^{reg}@>\si>>\cb_u^{reg}$.}
\nl
If $P,P'\in\cp_u^{reg}$ and $\ph(P)=\ph(P')=B$ then $P,P'$ are parabolic subgroups in 
$\cp_{1245}$ containing $B$, hence $P=P'$. Assume now that $P\in\cp_u^{reg}$ and let 
$B=\ph(P)$. For $i\in\{1,2,4,5\}$ let $P_i$ be the parabolic subgroup in $\cp_i$ such that 
$B\sub P_i$. We have $P_i\sub P$ and $u\in U_B$. Since $\p_P(u)$ is regular unipotent in $\bP$ 
we have $\p_P(u)\n\p_P(U_{P_i})$ that is, $u\n U_{P_i}$; but this is the same as saying that 
$B\ne\cb_{u,i}$. Thus, $\ph(P)\in\cb_u^{reg}$. We see that $\ph$ restricts to an injective map 
$\ph_0:\cp_u^{reg}@>>>\cb_u^{reg}$. Now let $B\in\cb_u^{reg}$. Let $P$ be the unique parabolic 
subgroup in $\cp_{1245}$ such that $B\sub P$. We have $u\in P$. Again, for $i\in\{1,2,4,5\}$, 
let $P_i$ be the parabolic subgroup in $\cp_i$ such that $B\sub P_i$. We have $P_i\sub P$ and 
$u\in U_B$. Since $B\n\cb_{u,i}$, we have $u\n U_{P_i}$ hence $\p_P(u)\n\p_P(U_{P_i})$. It 
follows that $\p_P(u)$ is regular unipotent in $\bP$. Thus $\ph_0$ is surjective hence a 
bijection. We omit the proof of the fact that $\ph_0$ is an isomorphism.

\subhead 4.24\endsubhead
We show:

(a) {\it Let $\cl$ be a $0$-line in $\cb$ such that $\cl\cap\cb_u^{reg}\ne\emp$. If 
$\cl\sub\cb_u$, then $\cl\sub\cb_u^{reg}$. If $\cl\not\sub\cb_u$, then 
$\cl\cap\cb_u=\cl\cap\cb_u^{reg}$ is a single point.}
\nl
Assume first that $\cl\sub\cb_u$. Let $B\in\cl\cap\cb_u^{reg}$. If $\cl$ contains some 
$B'\ne B$ such that $B'\in\cb_{u,i}$ for some $i\in\{1,2,4,5\}$ then, applying 1.1(a) with $B$ 
replaced by $B'$ and $J=\{0,i\}$ we see that $\cb^P\sub\cb_u$ where $P\in\cp_{0,i}$ contains 
$B'$. We have $B\in\cb^P$ and the line of type $i$ through $B$ is contained in $\cb^P$ hence in
$\cb_u$, so that $B\in\cb_{u,i}$; this contradicts $B\in\cb_u^{reg}$. We see that 
$\cl\sub\cb_u^{reg}$.

Assume next that $\cl\not\sub\cb_u$. Clearly, the intersection of any $0$-line with $\cb_u$ is 
either empty, or $\cl$, or a point. In our case $\cl\cap\cb_u$ is not $\cl$ and is not empty 
hence it is a point. Hence $\cl\sub\cb_u^{reg}$ is either empty or a point; by assumption it is
nonempty hence it is a point. This proves (a).

\subhead 4.25\endsubhead
From 4.24(a) we see that we have a partition $\cb_u^{reg}={}'\cb_u^{reg}\cup{}''\cb_u^{reg}$ 
where ${}''\cb_u^{reg}$ is the union of the $0$-lines contained in $\cb_u^{reg}$ and 
${}'\cb_u^{reg}$ is the set of all $B\in\cb_u^{reg}$ such that the $0$-line through $B$ 
intersects $\cb_u$ in exactly one point, $B$. Note that ${}'\cb_u^{reg}=\cb_u^\di$. Hence, 
using 4.22(b), we have

(a) {\it $\c({}'\cb_u^{reg})\in\{0,1,2\}$.}
\nl
For any $P\in\cp_{1245}$ we denote by $Q_P$ the unique parabolic subgroup in $\cp_{01245}$
such that $P\sub Q_P$; note that $(\ov{Q_P})_{ad}=\bP_{ad}\T H_P$ (canonically) where 
$H_P\cong PGL_2(\kk)$ (we use that $s_0$ commutes with $W_{1245}$). Let ${}''\cp_u^{reg}$ 
(resp. ${}'\cp_u^{reg}$) be the set of all $P\in\cb_u^{reg}$ such that the image of $u\in Q_P$ 
under the obvious composition $Q_P@>>>(\ov{Q_P})_{ad}@>>>H_P$ is $1\in H_P$ (resp. a regular 
unipotent element of $H_P$).

From the definitions we see that under the isomorphism $\ph_0:\cp_u^{reg}@>\si>>\cb_u^{reg}$ in
4.23(a), the subset ${}''\cp_u^{reg}$ of $\cp_u^{reg}$ corresponds to the subset
${}''\cb_u^{reg}$ of $\cb_u^{reg}$; hence ${}'\cp_u^{reg}=\cp_u^{reg}-{}''\cp_u^{reg}$ 
corresponds under $\ph_0$ to ${}'\cb_u^{reg}=\cb_u^{reg}-{}''\cb_u^{reg}$. Using this and (a) 
we deduce

(b) {\it $\c({}'\cp_u^{reg})\in\{0,1,2\}$.}

\subhead 4.26\endsubhead
We now assume that $G$ is simply connected and $p\ne3$. Let $J=\{1,2,4,5\}\sub I$. We fix
$P_J\in\cp_J$. Let $C$ be the regular unipotent conjugacy class of $\bP_J$. Let $\cs_0$ be an 
irreducible cuspidal $\bP_J$-equivariant local system on $C$. Up to isomorphism there are two 
such local systems, one for each nontrivial character of the group (cyclic of order $3$) of 
connected components of the centralizer in $\bP_J$ of an element in $C$. We shall use the 
notation and results of 2.1, 2.2 for this $G,J,P_J,\cs_0$.

In our case $\cw$ in 2.1 is a Weyl group of type $G_2$ with simple reflections $\s_0,\s_3$.
Now $\XX_u$ and the local system $\hcs$ on it are defined as in 2.1 (with $g=u$).
Recall from 2.1 that $H^j_c(\XX_u,\hcs)$ is naturally a $\cw$-module for $j\in\ZZ$.
We have the following result.
$$\sum_j(-1)^j\tr(\s_0,H^j_c(\XX_u,\hcs))\in\{0,1,2\}.\tag a$$ 
Let $J'=J\cup\{0\}$ and let $M_{J',g,r},\p'',\dcs_r$ be as in 2.1, 2.2 with $g=u$, $i=0$.
We set $\p''{}\i(M_{J',u,r})=\cm_r$.
From 2.2(b) we see that the left hand side of (a) is equal to $\c(\cm_r,\dcs_r)$
From the definitions we see that there exists an unramified principal covering
$\ps:\hat\cm_r@>>>\cm_r$ with group $\ZZ/3$ (with generator $\k$) and 
$\th\in\bbq-\{1\}$ with $\th^3=1$ such that the stalk of $\dcs_r$ at any $x\in\cm_r$
is equal to the vector space of functions $f:\ps\i(x)@>>>\bbq$ such that $f(\k\tx)=\th f(\tx)$ 
for any $\tx\in\ps\i(x)$. It follows that 
$$\c(\cm_r,\dcs_r)=\sum_{h=0}^2\sum_j(-1)^j\tr((\k^h)^*,H^j_c(\hat\cm_r,\bbq))\th^h/3,$$
$$\c(\cm_r)=\sum_{h=0}^2\sum_j(-1)^j\tr((\k^h)^*,H^j_c(\hat\cm_r,\bbq))/3,$$
where $(\k^h)^*$ is induced by the action of $\k^h$ on $\hat\cm_r$. If $h\in\{1,2\}$ then 
$\k^h:\hat\cm_r@>>>\hat\cm_r$ has no fixed points and it has order $3$. Since $p\ne3$ we can 
apply the fixed point formula in \cite{\DL, 3.2} to see that
$\sum_j(-1)^j\tr((\k^h)^*,H^j_c(\hat\cm_r,\bbq))=0$ for $h\in\{1,2\}$. It follows that
$$\c(\cm_r,\dcs_r)=\sum_j(-1)^j\dim H^j_c(\hat\cm_r,\bbq))/3,$$
$$\c(\cm_r)=\sum_j(-1)^j\dim H^j_c(\hat\cm_r,\bbq))/3,$$
so that 
$$\c(\cm_r,\dcs_r)=\c(\cm_r).$$
Thus, to prove (a), it is enough to show that
$$\c(\cm_r)\in\{0,1,2\}.$$ 
From the definitions we see that the map $\cm_r@>>>{}'\cp_u^{reg}$, 
$(xP_J,y(xU_{J'}x\i))\m xP_Jx\i$, is a isomorphism. Hence it is enough to show that
$\c({}'\cp_u^{reg})\in\{0,1,2\}$. But this is known from 4.25(b). This completes the proof of
(a).

\head 5. The main result\endhead
\subhead 5.1\endsubhead
In this section we assume that $G$ is simply connected of type $E_6$ and that $p\ne3$. We
write 
$I=\{0,1,2,3,4,5\}$ as in 4.1. Let $J=\{1,2,4,5\}\sub I$ and let $P_J\in\cp_J$. Let $C$ be the 
regular unipotent class in $\bP_J$. There are exactly two irreducible cuspidal $\bP_J$-equivariant
local systems on $C$; we fix one of them, say $\cs_0$. Let $\cw=N_WW_J/W_J$ where $N_WW_J$ is the 
normalizer of $W_J$ in $W$. This is a finite Coxeter group of type $G_2$ with simple reflections 
$\s_0,\s_3$ defined as in 4.1. Now the block defined by $(P_J,C,\cs_0)$ consists of six pairs
$(C_h,\ce_h)$, $h\in H:=\{0,1,3,4,9,12\}$, where $C_h$ is a unipotent class of $G$ such that, for 
$u_h\in C_h$ we have $d_{u_h}=h$ and $\ce_h$ is an irreducible $G$-equivariant local system on 
$C_h$, uniquely determined up to isomorphism by $h$. In  \cite{\SPAA}, 
$C_0,C_1,C_3,C_4,C_9,C_{12}$ are denoted by $E_6, E_6(a_1),A_5A_1,A_5,2A_2A_1,2A_2$ respectively. We 
shall assume that $u\in G$ (see 1.1) belongs to $C_3$.

Let $\Irr\cw$ be the set of irreducible representations (up to isomorphism) of $\cw$. We have
$\Irr\cw=\{1,s,\e,\e',\r,\r'\}$ where
$1,s,\e,\e'$ are one-dimensional, $\r,\r'$ are two-dimensional, $1$ is the unit representation,
and we have $\r\ot\r=\r'\ot\r'=1+\r'+s$, $\r\ot\r'=\e+\e'+\r$; these properties distinguish
$\r$ from $\r'$ and $s$ from $1,\e,\e'$ but not $\e$ from $\e'$. 
We can distinguish $\e$ from $\e'$ by the following requirements:
$\s_0$ acts as $-1$ on $\e$ and as $1$ on $\e'$; $\s_3$ acts as $-1$ on $\e'$ and as $1$ on $\e$. 
We have also
$\e\ot\e=\e'\ot\e'=s\ot s=1$, $\e\ot\e'=s$, $\e\ot s=\e'$, $\e'\ot s=\e$,
$\r\ot s=\r$, $\r'\ot s=\r'$, $\r\ot\e=\r\ot\e'=\r'$, $\r'\ot\e=\r'\ot\e'=\r$.

The generalized Springer correspondence for our block provides a bijection 
$\Irr\cw\lra\{(C_h,\ce_h);h\in H\}$. For $h\in H$ let $E_h\in\Irr\cw$ be corresponding to
$(C_h,\ce_h)$ under this bijection. According to \cite{\SPAA}, either (a) or (b) below holds:

(a) $E_0=1$, $E_1=\e$, $E_3=\r$, $E_4=\r'$, $E_9=\e'$, $E_{12}=s$;

(b) $E_0=1$, $E_1=\e$, $E_3=\r'$, $E_4=\r$, $E_9=\e'$, $E_{12}=s$.

\subhead 5.2\endsubhead
Let $q$ be an indeterminate. For any representation $E$ of $\cw$ we set
$$\G_E=(1/12)\sum_{w\in\cw}(q^2-1)(q^6-1)\tr(w,E)/\det(q-w)\in\NN[q].$$

We have
$$\G_1=q^6,\G_\e=\G_{e'}=q^3,\G_\r=q+q^5,\G_{\r'}=q^2+q^4,\G_s=1.$$. It follows that the matrix 
$(\G_{E_h\ot E_{h'}})_{(h,h')\in H\T H}$ is equal to
$$\matrix
q^6 &    q^3 &          q+q^5 &      q^2+q^4 &    q^3&    1\\
q^3 &    q^6 &         q^2+q^4&      q+q^5 &      1 &     q^3\\
q+q^5&   q^2+q^4&   1+q^2+q^4+q^6&  q+2q^3+q^5&  q^2+q^4&  q+q^5\\
q^2+q^4&  q+q^5&    q+2q^3+q^5&    1+q^2+q^4+q^6& q+q^5&  q^2+q^4\\
q^3&      1&          q^2+q^4&       q+q^5&       q^6&     q^3\\
1 &      q^3&          q+q^5 &       q^2+q^4&     q^3&     q^6\endmatrix$$
if 5.1(a) holds and to
$$\matrix
q^6 &    q^3 &          q^2+q^4 &      q+q^5 &    q^3 &   1\\
q^3 &    q^6 &         q+q^5  &      q^2+q^4 &      1 &     q^3\\
q^2+q^4&  q+q^5&   1+q^2+q^4+q^6&  q+2q^3+q^5&  q+q^5&  q^2+q^4\\
q+q^5&  q^2+q^4&    q+2q^3+q^5&    1+q^2+q^4+q^6& q^2+q^4&  q+q^5\\
q^3 &     1&          q+q^5&       q^2+q^4 &      q^6 &    q^3\\
1   &    q^3&          q^2+q^4&        q+q^5&     q^3&     q^6\endmatrix$$
if 5.1(b) holds. It follows that the matrix $M:=(q^{-h-h'}\G_{E_h\ot E_{h'}})_{(h,h')\in H\T H}$ 
is equal to 
$$\matrix
q^6  &        q^2  &         q^{-2}+q^2  &     q^{-2}+1 &        q^{-6} &   q^{-12}\\
q^2 &         q^4  &        q^{-2}+1    &     q^{-4}+1  &        q^{-10}&    q^{-10}\\
q^{-2}+q^2 &q^{-2}+1& q^{-6}+q^{-4}+q^{-2}+1& q^{-6}+2q^{-4}+q^{-2}& 
        q^{-10}+q^{-8}&  q^{-14}+q^{-10}\\
q^{-2}+1& q^{-4}+1& q^{-6}+2q^{-4}+q^{-2}& q^{-8}+q^{-6}+q^{-4}+q^{-2} &
          q^{-12}+q^{-8}&  q^{-14}+q^{-12}\\
q^{-6}& q^{-10}&          q^{-10}+q^{-8}&   q^{-12}+q^{-8}& q^{-12}& q^{-18}\\
q^{-12}& q^{-10}& q^{-14}+q^{-10}& q^{-14}+q^{-12}& q^{-18}&     q^{-18}\endmatrix$$
if 5.1(a) holds and to
$$\matrix
q^6 &    q^2&           q^{-1}+q&     q^{-3}+q&     q^{-6}&    q^{-12}\\
q^2 & q^4 &    q^{-3}+q&        q^{-3}+q^{-1}&    q^{-10}&   q^{-10}\\
q^{-1}+q &q^{-3}+q& q^{-6}+q^{-4}+q^{-2}+1& q^{-6}+2q^{-4}+q^{-2}& q^{-11}+q^{-7} &
q^{-13}+q^{-11}\\
q^{-3}+q &q^{-3}+q^{-1}& q^{-6}+2q^{-4}+q^{-2}& q^{-8}+q^{-6}+q^{-4}+q^{-2} &
q^{-11}+q^{-9}&  q^{-15}+q^{-11}\\ 
q^{-6}&     q^{-10}&   q^{-11}+q^{-7}&   q^{-11}+q^{-9}&   q^{-12}&     q^{-18}\\
q^{-12}& q^{-10}&   q^{-13}+q^{-11}&     q^{-15}+q^{-11}&  q^{-18}&     q^{-18}\endmatrix$$
if 5.1(b) holds.

\subhead 5.3\endsubhead
We can write uniquely $M={}^t\Pi A\Pi$ where $\Pi,A$ are matrices indexed by $H\T H$ such that
$A$ is diagonal and ${}^t\Pi$ is upper triangular with $1$ on diagonal. More precisely $A$ is equal
to
$$\matrix
q^{-2}(q^2-1)(q^6-1)&       0&              0&           0&         0&        0\\
0     &              q^{-4}(q^2-1)(q^6-1)&  0&           0&         0&        0\\
0     &                0&          q^{-8}(q^2-1)(q^6-1)& 0 &        0&        0\\
0   &                  0&                   0&      q^{-8}(q^6-1)&   0&        0\\
0   &                  0&                   0&           0&    q^{-18}(q^6-1)& 0 \\
0   &                  0&                   0&           0&          0 &     q^{-18}
\endmatrix$$
both in case 5.1(a) and in case 5.1(b); moreover ${}^t\Pi$ is equal to
$$\matrix
1  & 0 &   q^2&    q^2&    q^6&     q^6\\
0 &  1 &   0&      q^2&    0 &      q^8\\
0  & 0 &   1&      1 &    q^4&      q^4+q^8\\
0  & 0 &   0 &     1 &    q^4&      q^4+q^6\\
0 &  0 &   0 &     0 &     1&       1\\      
0 &  0 &   0 &     0 &     0&       1\endmatrix$$      
if 5.1(a) holds and to
$$\matrix
1 &  0 &   0 &     q^3&    q^6&     q^6\\
0 &  1 &   q &     q &    0 &      q^8\\
0 &  0 &   1 &     1 &    q^5 &     q^5+q^7\\
0 &  0 &   0 &     1 &    q^3 &     q^3+q^7\\
0 &  0 &   0 &     0 &     1  &     1\\
0 &  0 &   0 &     0 &     0 &      1\endmatrix$$
if 5.1(b) holds.

Let $h\le h'$ in $H$; let $\Pi_{h',h}$ be the $(h',h)$-entry of $\Pi$ or the $(h,h')$-entry of 
${}^t\Pi$. We write $\Pi_{h',h}=\sum_{k\in\NN}\Pi^k_{h',h}q^k\in\NN[q]$ where $\Pi^k_{h',h}\in\NN$.

Note that $C_{h'}$ is contained in the closure $\bC_h$ of $C_h$. For $k'\in\ZZ$ let 
$$\ch^{k'}IC(\bC_h,\ce_h)|_{C_{h'}}$$ 
be the restriction to $C_{h'}$ of the $k'$-th cohomology sheaf
of the intersection cohomology complex $IC(\bC_{h'},\ce_{h'})$; this is a $G$-equivariant local 
system on $C_h$. 

In the remainder of this subsection we assume that  $p\n\{2,3\}$. 
From \cite{\CSV, 24.8} we have that 
$$\ch^{k'}IC(\bC_h,\ce_h)|_{C_{h'}}=0$$ 

for $k'$ odd and that for $k\in\ZZ$, 

$\ch^{2k}IC(\bC_h,\ce_h)|_{C_{h'}}$ is isomorphic to a direct sum of $\Pi^k_{h',h}$ 
copies of $\ce_{h'}$. 

\subhead 5.4\endsubhead
In this subsection we assume that $p\n\{2,3\}$. 

Recall from 2.1, 2.2 that
$$\XX_u=\{xP_J\in G/P_J;x\i gx\in P_J,\p_J(x\i ux)\in C\}$$
carries a local system $\hcs$ induced from $\cs_0$ and that for any $j\in\ZZ$, $\cw$ acts 
naturally on $H^j_c(\XX_u,\hcs)$. (Note that $H^j_c(\XX_u,\hcs)=0$ for odd $j$.) Let 
$m_{j,h}$ be the number of times the irreducible $\cw$-module $E_h$ appears in the $\cw$-module
$H^{2j}_c(\XX_u,\hcs)$. Using \cite{\ICC, 6.5}, we see that for $h\le3$ we have
$\sum_jm_{j,h}q^j=\Pi_{3,h}q^h\in\ZZ[q]$. Using the tables in 5.3 we deduce that, as a 
$\cw$-module, $H^6_c(\XX_u,\hcs)$ is $\r$ if 5.1(a) holds and is $\r'$ if 5.1(b) holds;
$H^4_c(\XX_u,\hcs)$ is $1$ if 5.1(a) holds and is $\e$ if 5.1(b) holds; $H^j_c(\XX_u,\hcs)$ is $0$
if $j\ne\{4,6\}$.
It follows that 

(a) $\sum_j(-1)^j\tr(\s_0,H^j_c(\XX_u,\hcs))$ is equal to
$\tr(\s_0,\r)+1=1$ if 5.1(a) holds and to $\tr(\s_0,\r')-1=-1$ if 5.1(b) holds.
\nl
From 4.26 we see that $\sum_j(-1)^j\tr(\s_0,H^j_c(\XX_u,\hcs))=\c({}'\cp_u^{reg})$
where ${}'\cp_u^{reg}$ is as in 4.25. Thus, $\c({}'\cp_u^{reg})$ is equal to $1$ if
5.1(a) holds and to $-1$ if 5.1(b) holds.
By 4.25(b) we have $\c({}'\cp_u^{reg})\in\{0,1,2\}$. We deduce the following result in which
we assume that $p\n\{2,3\}$.

\proclaim{Theorem 5.5}(a) Statement 5.1(a) holds.

(b) We have $\c({}'\cp_u^{reg})=1$.
\endproclaim

\subhead 5.6\endsubhead
One can show that the statement of 5.5 remains true when $p=2$. In this case the references to 4.26 
and 4.25(b) can still be used. Although the statements at the end of 5.3 are not known in this case,
the statement 5.4(a) remains valid in this case.

\head 6. Final comments\endhead
\subhead 6.1\endsubhead
In this and the next subsection 
we assume that $G$ is simply connected/adjoint of type $E_8$ and that $p=3$. We
write $I=\{0,1,2,3,4,5,6,7\}$ where the numbering is chosen so that 
$s_1s_2,s_2s_3,s_3s_4,s_4s_5,s_5s_6,s_6s_7,s_0s_3$ have order $3$. 
Let $J=\{0,1,2,3,4,5\}\sub I$ and let $P_J\in\cp_J$. Let $C$ be the 
regular unipotent class in $\bP_J$. There are exactly two irreducible cuspidal $\bP_J$-equivariant
local systems on $C$; we fix one of them, say $\cs_0$. Let $\cw=N_WW_J/W_J$ where $N_WW_J$ is the 
normalizer of $W_J$ in $W$. This is a finite Coxeter group of type $G_2$ with simple reflections 
$\s_6,\s_7$ defined as in 4.1. Now the block defined by $(P_J,C,\cs_0)$ consists of six pairs
$(C_h,\ce_h)$, $h\in H:=\{0,1,3,4,9,12\}$, where $C_h$ is a unipotent class of $G$ such that, for 
$u_h\in C_h$ we have $d_{u_h}=h$ and $\ce_h$ is an irreducible $G$-equivariant local system on 
$C_h$, uniquely determined up to isomorphism by $h$. In  \cite{\SPAA}, 
$C_0,C_1,C_3,C_4,C_9,C_{12}$ are denoted by 
$E_8,E_8(a_1),E_7A_1,E_7,E_6A_1,E_6$ respectively. We shall assume that $u\in G$ (see 1.1) belongs to 
$C_3$. We identify the present $\cw$ with the group denoted by $\cw$ in 5.1 by requiring that
$\s_6,\s_7$ correspond respectively to $\s_3,\s_0$ in 5.1. Note that $s_7$ commutes with $W_J$.

The generalized Springer correspondence for our block provides a bijection 
$\Irr\cw\lra\{(C_h,\ce_h);h\in H\}$. For $h\in H$ let $E_h\in\Irr\cw$ be corresponding to
$(C_h,\ce_h)$ under this bijection. According to \cite{\SPAA}, either (a) or (b) below holds:

(a) $E_0=1$, $E_1=\e$, $E_3=\r$, $E_4=\r'$, $E_9=\e'$, $E_{12}=s$;

(b) $E_0=1$, $E_1=\e$, $E_3=\r'$, $E_4=\r$, $E_9=\e'$, $E_{12}=s$.
\nl
(Compare with 5.1(a),(b).)

\proclaim{Conjecture 6.2} Statement 6.1(a) holds.
\endproclaim
Recall from 2.1, 2.2 that
$$\XX_u=\{xP_J\in G/P_J;x\i gx\in P_J,\p_J(x\i ux)\in C\}$$
carries a local system $\hcs$ induced from $\cs_0$ and that for any $j\in\ZZ$, $\cw$ acts 
naturally on $H^j_c(\XX_u,\hcs)$. From \cite{\CSV} one can deduce, using computations like those in
5.2, that 

(a) $\sum_j(-1)^j\tr(\s_7,H^j_c(\XX_u,\hcs))$ is equal to
$\tr(\s_7,\r)+1=1$ if 6.1(a) holds and to $\tr(\s_7,\r')-1=-1$ if 6.1(b) holds.
\nl
Hence to prove the conjecture it is enough to show that 

(b) $\sum_j(-1)^j\tr(\s_7,H^j_c(\XX_u,\hcs))\in\NN$.
\nl
By arguments similar to those in 4.26 we see that the sum in (b) is equal to the Euler characteristic of
$$\cb_u-(\cb_{u,0}\cup\cb_{u,1}\cup\cb_{u,2}\cup\cb_{u,3}\cup\cb_{u,4}\cup\cb_{u,5}\cup\cb_{u,7})$$
with coefficient in a local system of rank $1$ defined by $\cs_0$. Hence it is enough to show that this
Euler characteristic is in $\NN$. It is not clear how to do this. One difficulty is that the argument 
in 4.26 (based on \cite{\DL}) is not applicable 
hence in the Euler characteristic above the local system cannot be replaced by $\bbq$.

\enddocument